%% file: draft_5.4_arxiv.tex
\documentclass{amsart}
\usepackage[letterpaper, portrait, margin=1in]{geometry}
\usepackage{style}

\usepackage[table]{xcolor}
\usepackage{colortbl}
\definecolor{c0}{HTML}{000000}
\definecolor{c1}{HTML}{9AB8FD}
\definecolor{c2}{HTML}{E361A0}
\definecolor{c3}{HTML}{FD8F60}
\definecolor{c4}{HTML}{998BF2}
\definecolor{c5}{HTML}{1CCA93}
\definecolor{c6}{HTML}{FEC161}
\definecolor{c7}{HTML}{D40000}

\input{macros}

\title[Special points on intersections of hypersurfaces]{Special points on intersections of hypersurfaces}

\author{Claudio G\'{o}mez-Gonz\'{a}les}
\address{Department of Mathematics \& Statistics, Carleton College}
\email{cgonzales@carleton.edu}

\begin{document}

\begin{abstract}
We establish lower bounds on the ambient dimension for an intersection of hypersurfaces to have a dense collection of ``level $\ell$'' points, in the sense introduced by Arnold--Shimura, given as a polynomial in the numbers of hypersurfaces of each degree. Our method builds upon \cite{GGW2025} to include other classes of accessory irrationality, towards the problem of understanding the arithmetic of ``special points'' as introduced in \cite{GGSW2024}. We deduce improved upper bounds on resolvent degree $\RD(n)$ and $\RD(G)$ for the sporadic groups as part of outlining frameworks for incorporating future advances in the theory.
\end{abstract}

\maketitle

\section{Introduction}\label{sec:introduction}

Resolvent degree is a measure of complexity motivated by classical phenomena in algebraic geometry which naturally generalizes the problem of solving algebraic equations in radicals. In its most basic form, formalized independently by Brauer \cite{Brauer1975} and Arnold--Shimura \cite{AS1976}, $\RD(n)$ asks for the minimal $r$ such that the generic degree $n$ polynomial can be solved in terms of its coefficients using algebraic functions of at most $r$ variables. The quadratic, cubic, and quartic formulas show that the roots of degree $\leq 4$ polynomials can be expressed via radicals---meaning they can be obtained using algebraic functions of a single variable---and give $\RD(2) = \RD(3) = \RD(4) = 1$. 

For polynomials of degree $5$, the Abel--Ruffini theorem states that a solution in radicals is impossible. However, Bring \cite{Bring1786} demonstrated that any quintic can be simplified (via radicals) to the form $x^5 + ax + 1$, showing that $\RD(5) = 1$. Hilbert's 13th Problem asks whether the bound $\RD(7) \leq 3$ is sharp, and his Sextic and Octic Conjectures ask the same question of the bounds $\RD(6) \leq 2$ and $\RD(8) \leq 4$, respectively. We emphasize that whether $\RD(n) > 1$ for any $n$ is an open problem. 

Historically, many upper bounds on $\RD(n)$ have been phrased in terms of values at which the number of algebraically independent coefficients in a generic polynomial can be systematically reduced (\cite{Tschirnhaus1683}). For example, Hamilton's \cite{Hamilton1836} extension of Bring's reduction of the quintic established that
\[ 
    \RD(5) = 1, \qquad \RD(6) \leq 2, \qquad \RD(7) \leq 3, \qquad \RD(8) \leq 4, 
\]
and so on, i.e., that at least $4$ parameters can be simplified from the generic degree $n$ polynomial whenever $n \geq 5$. 
Following \cite{Dixmier1993}, we formally introduce $\mch: \Nb_{\geq 1} \to \Nb_{\geq 1}$ as a function to capture this phenomenon:
\[ 
    \mch(r) \coloneqq \min \{ n \in \Nb_{\geq 1}: n-\RD(n) \geq r \}. 
\]

Establishing upper bounds on $\mch(r)$ has been implicit in much of the literature surrounding resolvent degree. 
From this perspective, Bring showed $\mch(4) = 5$, and Hilbert's sketch \cite{Hilbert1927}\footnote{The first proof of this fact is due to Wiman \cite{Wiman1927, Sutherland2021Translation}; a proof also appears in the appendix of \cite{Dixmier1993}.} of a 4-variable formula for the general nonic via the 27 lines on a cubic surface gives $\mch(5) \leq 9$. 
Among other things, Hamilton established that $\mch(r)$ is well-defined for all $r$. 
Modern interest begins with Brauer's bound \cite{Brauer1975}:
\begin{equation}\label{eqn:brauer}
    \mch(r) \leq (r-1)! + 1,
\end{equation}
for any $r \geq 4$, and was renewed more recently by Wolfson, who advanced Hilbert's technique into a more general geometric method \cite{Wolfson2020}. The best bounds prior to the present paper appear in \cite{Sutherland2022, HS2023}. Our first results improve upon these upper bounds for several values of $r$.

\begin{theorem}[Appears as Theorem \ref{cor:new_bounds_on_H_body}]\label{thm:new_bounds_on_H} With notation as above, the following bounds hold:
\begin{align*}
\mch(7) &\leq 74, & \mch(8) &\leq 211, & \mch(11) &\leq 59050, \\
\mch(12) &\leq 332641, & \mch(13) &\leq 3991681, & \mch(20) &\leq 227214539745187, \\
\mch(21) &\leq 3379030566912001, & & \text{and} & \mch(22) &\leq 70959641905152001.
\end{align*}
\end{theorem}

While the true values of $\mch(r)$ remain unknown, these results represent a substantial relative reduction from the best bounds in the literature. For example, our bound on $\mch(7)$ reduces the prior estimate by nearly a third; for $r=21$, the bound is reduced to a sixth of its previous value. 

\medskip
These improved bounds are obtained by building upon the geometric methods developed in \cite{HS2023, GGW2025}. The first of these recovers Sylvester's obliteration algorithm \cite{Sylvester1887} for solving systems of equations via the classical theory of polar cones. The second establishes bounds on the ambient dimension $N$ required to guarantee that the set of points rational over the solvable closure, $X(K^{\Sol}) \subset X$, is dense when $X$ is an intersection of quadrics, cubics, and quartics. 
In this paper, we study a broader setting where quadrics, cubics, and quartics are replaced by hypersurfaces of arbitrary degrees, and the solvable closure $K^{\Sol}$ is replaced by the level $\ell$ closure $K^{(\ell)}$---for a formal definition of the level $\ell$ closure, see \cite{AS1976, Reichstein2025} or Section \ref{subsec:ed_rd_special}.
More generally, we refer to points rational over such closures as \emph{special points}. Establishing the existence of special points on certain intersections (as in \cite{Wolfson2020} and Section \ref{sec:bounds_on_RD}) corresponds to changes of variables that simplify the generic degree $n$ polynomial. The specific closure used dictates the permitted \emph{accessory irrationalities} to enact these transformations (e.g., in radicals or in algebraic functions of at most $\ell$ variables).

Our approach involves systematizing Sylvester's algorithm, which produces points on an intersection using linear subspaces contained in some subintersection and, in turn, determines linear subspaces via finding points on auxiliary varieties. This method allows us to formulate general statements regarding $j$-dimensional linear subspaces of intersections defined over level $\ell$ closures. Writing $\mcf(j,X)$ for the Fano variety of $j$-planes on $X$, we show:

\begin{theorem}[Appears as Theorem \ref{thm:polynomial_bounds_body}]\label{thm:polynomial_bounds} 
Fix $n \geq 2$, a field $k$ with $\characteristic(k) = 0$ or $\characteristic(k) > n$, and $\ell \geq \RD_k(n)$. Then there exists a polynomial $q \in \Qb[j,m_1,m_2,\dots,m_n]$ of total degree $2^{n-1}$ such that, for any field $K$ over $k$ and $X \subseteq \Pb^N_K$ given as an intersection of $m_i$ degree $i$ hypersurfaces (as $i$ ranges from $1$ to $n$), $\mcf(j,X)(K^{(\ell)}) \subseteq \mcf(j,X)$ is dense whenever
\[ N \geq q(j,m_1,m_2,\dots,m_n). \]
\end{theorem}

Given the general structure of this framework---finding points via linear subspaces, and linear subspaces via points---the upper bounds on the minimal ambient dimension $N$ established herein are inherently recursive. This recursion requires a sequence of choices governed by combinatorial objects called \emph{types}. A structured pathway through these choices that restricts $N$ efficiently, and consequently yields new bounds on resolvent degree, is formalized by our greedy obliteration heuristic (see Section \ref{sec:bounds_on_RD}).

\medskip
Resolvent degree as a measure of complexity can be applied to other algebraic objects. In particular, for a finite group $G$, the resolvent degree $\RD(G)$ measures the complexity of a general polynomial with Galois group $G$, recovering the classical setting of generic degree $n$ polynomials when $G = S_n$ (see \cite{FW2019, Reichstein2021, Sutherland2023} for more on these perspectives and accounts of recent progress). In \cite{GGSW2024}, building on techniques in \cite{Reichstein2025}, a strategy is established for bounding $\RD(G)$ in terms of the density of special points for certain $G$-varieties. 

Our main results regarding the density of level $\ell$ points on intersections of hypersurfaces permit a uniform approach for bounding $\RD(G)$ for finite simple groups, building on the methods of \cite{GGSW2024}, in terms of classical invariant theory. 
By this process, we sharpen upper bounds for the sporadic groups, including the Leech lattice group $\MCL$; the pariah $\RU$; the Monster subgroups $\HE, \FIWH, \FIWF, \B$; and the Monster itself $\MO$ (see \cite[Corollary~1.2 and Section~4]{GGSW2024} for previous estimates and additional context).

\begin{theorem}[Appears as Corollary \ref{cor:sporadic_bounds_body}]\label{thm:sporadic_bounds} 
The following bounds hold:
\[ 
\begin{aligned}
    \RD(\MCL) & \leq 18, \qquad & \qquad \RD(\RU) & \leq 25, \qquad & \qquad \RD(\HE) & \leq 47, \\
    \RD(\FIWH) & \leq 775, & \RD(\FIWF) & \leq 778, & \RD(\B) & \leq 4364, \\
    \text{ and } && \RD(\MO) & \leq 196872.
    \end{aligned} 
\]
\end{theorem}

Lastly, we return to the question of asymptotic behavior for $\mch(r)$. Because it remains an open question whether $\RD(n) > 1$ for any $n$, the true growth of $\mch(r)$ is entirely unknown. The methods outlined herein do not seem to improve the current best asymptotic upper bounds, which are given in terms of the existence of $j$-planes on intersections of hypersurfaces and obtained via the combinatorics of moduli space dimensions \cite[Section~3.3]{Sutherland2021}. 

Indeed, the sensitivity in the bounds we establish for the requisite ambient dimension $N$ (in order to guarantee dense level $\ell$ points) highlights why improving these asymptotic upper bounds is difficult. Our sharpest lower bounds on $N$ vary dramatically depending on the level $\ell$ of permitted accessory irrationalities. For small $\ell$, Theorem \ref{thm:polynomial_bounds} yields a bound on $N$ governed by a polynomial of degree $2^{n-1}$, where $n$ is the highest degree among the hypersurfaces. On the other hand, as $\ell$ grows, the required dimension drops, bounded from below by a polynomial of degree $n$. The sharp transition between these regimes is exactly the domain of interest to improving specific bounds on $\mch(r)$ (see Appendix \ref{appendix:computing_f} for more on these computational difficulties). We hope this paper encourages engagement with these subtleties, and with the arithmetic of special points more generally.

\medskip
Section \ref{sec:background} establishes the needed background on resolvent degree and accessory irrationalities; the combinatorics of the obliteration algorithm; and varieties of $j$-planes. Section \ref{sec:special_points_on_variety_of_j-planes} introduces a function $f^\mce$ which determines the number of variables needed to guarantee dense special points with respect to a specific class $\mce$ for intersections of a given type. It then establishes workhorse lemmas to approximate $f^\mce$ for the classes of interest (in analogy to \cite[Section~4]{GGW2025}) and proves Theorem \ref{thm:polynomial_bounds}. In addition to a brief survey of existing bounds on resolvent degree, Sections \ref{sec:bounds_on_RD} and \ref{sec:bounds_for_sporadics} establish procedures for applying this framework to sharpen a given bound on $\RD(n)$ and to obtain upper bounds on $\RD(G)$ via invariant theory, proving Theorems \ref{thm:new_bounds_on_H} and \ref{thm:sporadic_bounds}, respectively.
Lastly, Appendix \ref{appendix:computing_f} describes heuristics and obstacles for the computations concerning $f^\mce$. In particular, we compare the bounds in Theorem \ref{thm:polynomial_bounds} with the best estimates available to our methods in order to highlight the distinct regimes of behavior as $n$ and $\ell$ grow.

\subsection{Acknowledgments}
We are grateful to Alexander Sutherland and Benson Farb for helpful comments. We also thank Jesse Wolfson, who was closely involved with the ideas that led to this paper and provided thorough feedback on an earlier draft, but who declined to be listed as a coauthor. Special thanks are extended to Joshua Batson and Jesse Wolfson for pointing out non-local improvements to the obliteration strategy found using the AI tool Claude Fable 5, which also led to an improved bound on $\mch(7)$.
We lastly thank two anonymous reviewers for generous feedback, which substantially improved the exposition of this paper. The author was supported in part by the National Science Foundation Division of Mathematical Sciences, Grant No. DMS-2418943.

\section{Background}\label{sec:background}

Fix a ground field $k$ throughout. By a $k$-variety, we will mean a $k$-scheme of finite type (possibly reducible or non-reduced). While much of this paper relies on geometric constructions where scheme-theoretic structure is present---especially here and in Section \ref{sec:bounds_for_sporadics}---our primary results concern rational points over field extensions $K/k$. Because $X(K) = X^{\red}(K)$ as sets for any $k$-variety $X$, the non-reduced structure does not restrict the existence of rational points. Consequently, when evaluating rational points, we routinely identify schemes with their underlying reduced closed subschemes. 

Finally, while our eventual numerical results in Sections \ref{sec:bounds_on_RD} and \ref{sec:bounds_for_sporadics} specialize to the classical case $k = \Cb$, the foundational framework holds over any base field $k$ satisfying the characteristic conditions specified in each statement.

\subsection{\texorpdfstring{$\ed$}{ed}, \texorpdfstring{$\RD$}{RD}, and special points}\label{subsec:ed_rd_special}

Here we recall the fundamentals of essential dimension and resolvent degree, as well as classes of accessory irrationalities $\mce$. While the theory can be equivalently formulated in the language of branched covers and field extensions, we prefer the latter for the sake of compatibility with the recent literature upon which we rely (e.g., \cite{Reichstein2025, GGW2025}). For more on the equivalence of definitions, see \cite[Section~2]{GGSW2024}, whose treatment we follow.

Motivated by the classical problem of solving polynomials via algebraic functions of few variables, we first define resolvent degree for field extensions and finite groups.

\begin{definition}[Essential dimension and resolvent degree] Let $L/K$ be a finite extension.
\begin{enumerate}
    \item The \emph{essential dimension of $L/K$} is given by
    \[
        \ed_k(L/K) \coloneqq \min_{E} \left\{ \trdeg_k(E) \mid L/K \text{ is defined over } E \right\}.
    \]
    \item The \emph{resolvent degree of $L/K$}, written $\RD_k(L/K)$, is the least $d$ such that there is a tower of finite extensions,
    \[
        K = K_0 \into K_1 \into \cdots \into K_r,
    \]
    and an embedding $L \into K_r$ over $K$, such that $\ed_k(K_{i+1}/K_i) \leq d$ for all $i$. 
\end{enumerate}
\end{definition}

\begin{remark}
Since $L/K$ is defined over $K$, the trivial tower $K = K_0$ and $L = K_1$ gives
\begin{equation}\label{eqn:RD_ed_trdeg_bound}
    \RD_k(L/K) \leq \ed_k(L/K) \leq \trdeg_k(K).
\end{equation}
\end{remark}

To capture the intrinsic complexity of a finite group $G$ over a ground field $k$, we define its resolvent degree as the supremum over all possible $G$-Galois extensions.

\begin{definition}[Resolvent degree for finite groups]
The \emph{resolvent degree of $G$} is given by
\[ 
    \RD_k(G) \coloneqq \sup_{L/K} \left\{ \RD_k(L/K) \mid K\ \text{is a $k$-field and}\ L/K\ \text{is $G$-Galois} \right\}.
\]
For any $n \in \Nb_{\geq 1}$, we set
\[ 
    \RD_k(n) \coloneqq \RD_k(S_n)
\]
and simply write $\RD(n) \coloneqq \RD_\Cb(n)$, which corresponds to the classical notion of resolvent degree described in the introduction and used throughout the paper.
\end{definition}

Note that the supremum is well-defined, since such an extension $L/K$ always exists. Indeed, if $L$ is the field of rational functions on the group algebra $k[G]$, then the regular action of $G$ on $L$ is faithful, and by Artin's theorem $K = L^G$ satisfies $\Gal(L/K) \cong G$.

\medskip
Klein proved that \emph{accessory irrationalities}---such as the auxiliary quadratic used in his solution to the quintic---were necessary to write a formula for the quintic in $1$-variable functions \cite{Klein1884}. This idea was generalized in \cite{FKW2023} to the formal notion of \emph{classes of accessory irrationalities} $\mce$. In \cite{GGSW2024}, the authors developed a framework of ``special points''---points defined over closures with respect to a given class $\mce$---to establish arithmetic criteria used in bounding $\RD_k(G)$. In this paper, we are primarily interested in special points corresponding to the \emph{level $\ell$ closure}.

\begin{definition}\cite[Definition 6.1]{Reichstein2025}
Let $K$ be a field containing $k$, fix an algebraic closure $\bar{K}$, and let $\ell \geq 1$. The \emph{level $\ell$ closure} $K \into K^{(\ell)}$ is the compositum of all finite extensions $L/K$ such that $L \subseteq \bar{K}$ and $\RD_k(L/K) \leq \ell$. 
\end{definition}

To establish the framework we will use in Section \ref{sec:special_points_on_variety_of_j-planes}, we conclude with the general definition of accessory irrationality and accompanying closures.

\begin{definition}[Accessory irrationality]\cite[Definitions 2.5 and 2.7]{GGSW2024}\label{definition:accessory_irrationality}
Let $\Fieldsk$ be the category of fields over $k$ and let $\Fin: \Fieldsk \to \Cat$ be the functor sending $K$ to the category of finite, semisimple commutative $K$-algebras.\footnote{Recall by Wedderburn's theorem that any such algebra is a finite product of finite extensions of $K$.}
A \emph{class of accessory irrationalities} $\mce$ is a subfunctor $\mce \subset \Fin$ such that:
\begin{enumerate}
    \item For all $K$, $\mce(K)\subset \Fin(K)$ is a full subcategory.
    \item For all $K$, we have $K\in\mce(K)$.
    \item If $E,E'\in \mce(K)$, then $E\otimes_K E'\in\mce(K)$.
    \item If $L/K$ is finite, $E/L$ is finite, and $K\into L\into E$ is in $\mce(K)$, then $E\in \mce(L)$.
\end{enumerate}
Furthermore, we say that $\mce$ is
\begin{itemize}
    \item \emph{saturated} if for all finite extensions $L/K$ of $k$-fields, $K\into L\into E$ in $\mce(K)$ implies that $L\in \mce(K)$.
    \item \emph{closed under extensions} if for all finite extensions $L/K$ of $k$-fields, $L\in\mce(K)$ and $E\in\mce(L)$ together imply that $K\into L\into E$ is in $\mce(K)$.
\end{itemize}
\end{definition}

\begin{definition}[$\mce$-closure and special points]\cite[Definition 2.17]{GGSW2024}
Let $\mce$ be a saturated class of accessory irrationalities and suppose $K/k$ is a field extension with a fixed algebraic closure $\bar{K}$. The \emph{$\mce$-closure of $K$} is the compositum $K^\mce \coloneqq K(S_\mce) \into \bar{K}$, where
\[
    S_\mce \coloneqq \left\{ K \into L \into \bar{K} \mid L \in \mce(K) \right\}.
\]
Given a variety $X$ over $K$, we refer to points of the set $X(K^{\mce})$ as \emph{special points} of $X$ with respect to $\mce$.
\end{definition}

\begin{example}\label{example:e-closures} We describe several saturated classes of accessory irrationalities and their closures:
\begin{itemize}
    \item $\Sol$, the class of \emph{solvable accessory irrationalities}, sends $K$ to the full subcategory of finite semisimple commutative $K$-algebras which split as products of solvable extensions of $K$. This class is closed under extensions and $K \into K^{\Sol}$ is the maximal solvable closure.
    \item $\RD_k^{\le \ell}$, the class of \emph{level $\ell$ accessory irrationalities}, sends $K$ to the category of all finite semisimple commutative $K$-algebras $A$ such that $\RD_k(A/K)\le \ell$. The class $\RD_k^{\le \ell}$ is closed under extensions, as shown in \cite[Lemma 2.7]{FW2019}. The closure \[K \into K^{(\ell)} \coloneqq K^{\RD_k^{\le \ell}}\] appears in the literature at least as early as \cite{AS1976}. 
\end{itemize}
Note that, because radical extensions have $\RD=1$, we have $K^{\Sol} \into K^{(\ell)}$ for all $\ell \geq 1$. Similarly, for all $\ell \leq \ell'$, since $\RD_k^{\leq \ell}$ is a subfunctor of $\RD_k^{\leq \ell'}$ we have $K^{(\ell)} \into K^{(\ell')}$.
\end{example}

\subsection{Monoid of types}

\begin{definition}\label{definition:type}
A variety $X\subset \Pb^N_K$ is \emph{defined by equations of degrees $d_1,\ldots,d_r$} if it arises as a (scheme-theoretic) intersection of hypersurfaces $X=\bigcap_{j=1}^r \Vb(f_j)$ with $f_j \in K[x_0,\dots,x_N]_{d_j}$. Its \emph{type} is $\bm=(m_1,\ldots,m_n)$, where $m_i$ denotes the number of degree $i$ hypersurfaces in the intersection, i.e.,
\[
    m_i \coloneqq \left|\left\{f_j \mid \deg f_j=i\right\}\right|.
\]
For $X$ an intersection of hypersurfaces of type $\bm$, we also say $X^{\red}$ is of type $\bm$.
\end{definition}

\begin{remark}\label{remark:type_not_unique}
The type of $X$ is not determined by $X^{\red}$---the same reduced variety can appear as the reduction of schemes of different types. Further, operations involving polar cones (i.e., taking derivatives) may result in some of the defining polynomials of an intersection becoming identically zero. In defining type in terms of the graded pieces $K[x_0,\dots,x_N]_{d_j}$, we formally treat a zero polynomial arising in this way as a form of degree $d_j$.
\end{remark}

As in \cite[Section 2]{GGW2025}, we write $\Mt \coloneqq C_c(\Nb_{\geq 1}, \Nb)$ for the set of sequences of naturals with finite support, understood as a monoid with respect to entrywise addition whose combinatorics reflect the geometric constructions to follow. We will make use of the $\ell^1$ norm
\[ |\bm| = \sum_{i=1}^\infty m_i \]
and write $\ba \leq \bm$ for the entrywise partial ordering. We will also often suppress the rightmost tail of zeroes: if $\bm = (m_1,m_2,\dots) \in \Mt$ has $m_n$ as its last non-zero entry, we simply write
\[
    \bm = (m_1,\dots,m_n). 
\] 
We recall and extend some elementary notions here.

\begin{remark}\label{rem:passing_to_hyperplanes}
In this paper, we will frequently consider intersections of the form $(0,m_2,\dots,m_n)$. We omit linear forms because they do not complicate the problem of establishing density of $j$-planes, i.e., rather than considering an intersection of type $(m_1,\dots,m_n)$ in $\Pb^N$, one can study an intersection of type $(0,m_2,\dots,m_n)$ in $\Pb^{N-m_1}$.
\end{remark}

\begin{definition}\cite[Definition 2.7]{GGW2025}\label{def:j-endomorphism}
We define $(-^1): \Mt \to \Mt$ by
\[ 
    \bm^1\coloneqq(m_1,\dots,m_n)^1 \coloneqq \left(\sum_{i=1}^n m_i, \dots, \sum_{i=d}^n m_i, \dots, m_n \right) 
\]
and extend to a semigroup of maps $(-^j): \Mt \to \Mt$ by setting
\[ 
    \bm^{j+1} \coloneqq (\bm^j)^1 \quad \text{ and } \quad \bm^0 \coloneqq \bm.
\]
\end{definition}

\begin{lemma}\cite[Lemma 2.8]{GGW2025}\label{lemma:j-endomorphism}
Let $\bm \in \Mt$. Then:
\begin{enumerate}[label=(\alph*)]
\item For any $\ba \in \Mt$, we have $(\bm+\ba)^j=\bm^j+\ba^j$.
\item For any $j \in \Nb$, we have 
\[ \bm^j = \left( \sum_{i=1}^n \binom{j+i-2}{i-1} m_i, \dots, 
    \sum_{i=d}^n \binom{j+i-d-1}{i-d} m_i, \dots, m_n \right). \]
\end{enumerate}
\end{lemma}

\noindent In particular, it is clear from the above that $\ba \leq \bm$ means that $\ba^j \leq \bm^j$ for all $j \in \Nb$ and similarly that $\bm^j \leq \bm^{j'}$ whenever $j \leq j'$, which we will use freely.

We will need three additional combinatorial results in Section \ref{subsec:bounds_on_f} which allow us to control the endomorphisms $\bm \mapsto \bm^j$. Towards this end, despite thinking of $\Mt$ as an additive monoid, we will write $\bm - \ba$ whenever it is safe to do so---namely, when $\ba \leq \bm$---and note that the proof of Lemma \ref{lemma:j-endomorphism}(a) carries through verbatim when the $m_i$ are allowed to be negative.

\begin{lemma}\label{lemma:batching_obliterations}
Let $\bm, \ba, \bb \in \Mt$. Suppose that $\ba \leq \bm$ and $\bb \leq (\bm - \ba)^{|\ba|}$. Then we have
\[
    (\bm - (\ba + \bb))^{|\ba|+|\bb|} \leq \left( (\bm - \ba)^{|\ba|} - \bb \right)^{|\bb|}.
\]
\end{lemma}

\begin{proof}
    We expand using Lemma \ref{lemma:j-endomorphism}(a) and $\bm^{r+s} = (\bm^r)^{s}$. 
    The difference of the right-hand side minus the left-hand side is $\bb^{|\ba|+|\bb|} - \bb^{|\bb|}$, which is non-negative in every entry, and this completes the proof.
\end{proof}

For the last two results, we write $\be_d \in \Mt$ for the type equal to $1$ in degree $d$ and $0$ elsewhere.

\begin{lemma}\label{lemma:skipping_rightmost_combinatorics}
Let $\bm = (m_1,\dots,m_n) \in \Mt$, $1 \leq d < n$, and $j \in \Nb_{\geq 0}$. If $m_n, m_d \geq 1$, then 
\[
    \bm^j \leq (\bm - \be_d)^{j+1}.
\]
\end{lemma}

\begin{proof}
First, because $\bm$ and $\bm - \be_d$ are equal in entries above $d$, it is clear that the $i$th entry of $\bm^j$ will be less than or equal to that of $(\bm - \be_d )^{j+1}$ for all $i > d$. 

Next, consider the $d$th entry. For $\bm^j$, we have
\[
    m_d + \sum_{s=d+1}^{n} \binom{j+s-d-1}{s-d} m_s,
\]
while the $d$th entry of $(\bm - \be_d)^{j+1}$ is
\[
    (m_d - 1) + \sum_{s=d+1}^{n} \binom{j+s-d}{s-d} m_s.
\]
By Pascal's identity, the difference of these terms is:
\[
    -1 + \sum_{s=d+1}^{n-1} \binom{j+s-d-1}{s-d-1} m_s + \binom{j+n-d-1}{n-d-1} m_n.
\]
Since $m_n \geq 1$, the rightmost term is at least $1$; the intermediary terms are all non-negative, so the difference is strictly non-negative and the desired inequality holds at index $d$.

Next, consider $1 \leq i < d$. The $i$th entry of $(\bm - \be_d)^{j+1}$ is
\[
    m_i + \sum_{s=i+1}^{n} \binom{j+s-i}{s-i} m_s - \binom{j+d-i}{d-i}.
\]
Taking the difference with the $i$th entry of $\bm^j$ and applying Pascal's identity again gives
\[
    -\binom{j+d-i}{d-i} + \sum_{s=i+1}^{n-1} \binom{j+s-i-1}{s-i-1} m_s + \binom{j+n-i-1}{n-i-1} m_n.
\]
We must show this expression is non-negative; as before, because $m_n \geq 1$ and the intermediary terms are non-negative, this follows since $n-1 \geq d$:
\[ 
    \binom{j+n-i-1}{n-i-1} = \binom{j+n-i-1}{j} \geq \binom{j+d-i}{j} = \binom{j+d-i}{d-i}.
\]
\end{proof}

\begin{proposition}\label{prop:skipping_rightmost_combinatorics}
Let $\bm = (m_1,\dots,m_n), \ba = (a_1, \dots, a_{n-1}) \in \Mt$, and $j \in \Nb_{\geq 0}$. Suppose that $\ba \leq \bm$ and $j < m_n$. Then we have
\[
    (\bm - j \be_n )^j \leq (\bm - j \be_n - \ba)^{j+|\ba|}.
\]
\end{proposition}

\begin{proof}
    We proceed by induction on $|\ba|$, where the base case $\ba = 0$ is tautological.

    For the inductive step, suppose the inequality holds for all types of norm $r$, and consider a type $\ba \leq \bm$ with $|\ba| = r+1$. We write $\ba = \bb + \be_d$ for some $1 \leq d < n$ and $\bb$ of norm $r$. By the inductive hypothesis,
    \[
        (\bm - j \be_n)^j \leq (\bm - j \be_n - \bb)^{j+r}.
    \]
    Next, consider $\bm - j \be_n - \bb$, whose $n$th coordinate is $m_n - j \geq 1$. By Lemma \ref{lemma:skipping_rightmost_combinatorics}, we have
    \[
        (\bm - j \be_n - \bb)^{j+r} \leq (\bm - j \be_n - \bb - \be_d)^{j+r+1} = (\bm - j \be_n - \ba)^{j+r+1}.
    \]
    Taken together with the inductive hypothesis, we have shown the desired result.
\end{proof}

\subsection{Polar cones and Fano varieties of $j$-planes}

We next recall the fundamentals of polar cones. We refer the reader to \cite[Section 3]{GGW2025} in which the classical notion of polar cones (see \cite{Bertini1923} and \cite[Ch.~1, \S~Historical Notes]{Dolgachev2012}) and iterated constructions appearing in \cite{Sutherland2021, HS2023} are cast in terms of varieties of $j$-planes.

Firstly, recall that the \emph{Fano variety of $j$-planes on $X \subseteq \Pb^N$}, denoted $\mcf(j,X)$, is the subvariety of the Grassmannian $\grass(j,\Pb^N)$ of all $j$-planes contained in $X$, i.e.
\[
    \mcf(j,X) \coloneqq \{\Lambda\in \grass(j,\Pb^N) \mid \Lambda\subset X\}.
\]
To avoid confusion with the standard notion of ``Fano varieties''---i.e., complete varieties with ample anticanonical bundle---we refer to $\mcf(j,X)$ as simply the ``variety of $j$-planes on $X$'' henceforth. 

Given a fixed $\Lambda \in\mcf(j,X)$ and $i \geq j$, we will make use of subvarieties of $i$-planes containing $\Lambda$, written
\[
    \mcf_\Lambda(i,X) \coloneqq \{\Lambda' \in \mcf(i,X) \mid \Lambda \subseteq \Lambda' \}.
\]
Lastly, given $x_0,\ldots,x_j\in\Pb^N$, we write $\Lambda(x_0,\ldots,x_j)\subset \Pb^N$ for the plane they span.  

\medskip
Now let $X\subseteq \Pb^N$ be defined by equations of type $\bm$, and let $x_0\in X$. The polar cone $C^1(X;x_0)$ of a variety $X$ through $x_0 \in X$ is a subscheme of $X$ whose reduced subvariety is the union of all lines in $X$ through $x_0$, i.e., any $x_1 \in C^1(X;x_0) \setminus \{x_0\}$ defines a line $\Lambda(x_0,x_1) \subseteq X$. Moreover, if $X \subseteq \Pb^N$ has type $\bm$, then $C^1(X;x_0)$ is a variety of type $\bm^1$. As in Remark \ref{remark:type_not_unique}, passing to the polar cone may result in some of the new defining polynomials vanishing identically; we treat these zero polynomials as elements of the corresponding graded pieces $K[x_0,\dots,x_N]_{d_j}$ to formally determine the type $\bm^1$.

More generally, a \emph{$j$-polar point of $X$ defined with respect to $K/k$} is a tuple $(x_0,\dots,x_j)\in X^{j+1}(K)$ such that $\Lambda(x_0,\dots,x_j)$ is a $j$-plane contained in $X$. Given such a $j$-polar point, the $(j+1)$st polar cone of $X$ at $x_0,\dots,x_j$ is the variety 
\[ C^{j+1}(X;x_0,\dots,x_j) \coloneqq C^1(C^j(X;x_0,\dots,x_{j-1});x_j). \]
Iterated polar cones are nested by definition,
\[ C^{j+1}(X;x_0,\dots,x_j) \subseteq C^j(X;x_0,\dots,x_{j-1}) \subseteq \cdots \subseteq C^1(X;x_0) \subseteq X \subseteq \Pb^N, \]
and in general each $x_j \in C^j(X;x_0,\dots,x_{j-1}) \setminus \Lambda(x_0,\dots,x_{j-1})$ 
determines a $j$-plane $\Lambda(x_0,\dots,x_j) \subseteq X$.

\begin{theorem}[{\cite[Propositions~3.10 and 3.11]{GGW2025}}]\label{thm:non-empty_fano_dim}
Let $X$ be an intersection of type $\bm$ and let $(x_0,\dots,x_j)$ be a $j$-polar point of $X$ defined over $K/k$ with $\Lambda = \Lambda(x_0,\dots,x_j)$. 
Then $C^j(X;x_0,\dots,x_{j-1})^{\red}$ is the union of all $j$-planes in $X$ through $\Lambda(x_0,\dots,x_{j-1})$. Moreover,
\begin{enumerate}[label=(\alph*)]
    \item Writing $\bm^j$ as in Definition \ref{def:j-endomorphism}, $C^j(X;x_0,\dots,x_{j-1})$ has type $\bm^j$. In particular, $C^j(X;x_0,\dots,x_{j-1}) \setminus \Lambda(x_0,\dots,x_{j-1})$ is non-empty whenever $N \geq |\bm^j|+j$.
    \item For any fixed $(N-j)$-plane $\Lambda' \subset \Pb^N \setminus \Lambda(x_0,\dots,x_{j-1})$, the map
    \[
    \begin{aligned}
        C^j(X;x_0,\dots,x_{j-1}) \cap \Lambda' & \to \mcf_{\Lambda(x_0,\dots,x_{j-1})}(j,X) \\
        y & \mapsto \Lambda(x_0,\dots,x_j,y)
    \end{aligned}
    \]
    gives an isomorphism $C^j(X;x_0,\dots,x_{j-1})^{\red} \cap \Lambda' \xrightarrow{\cong} \mcf_{\Lambda(x_0,\dots,x_{j-1})}(j,X)^{\red}$.
\end{enumerate}
In particular, $\mcf_{\Lambda(x_0,\dots,x_{j-1})}(j,X)^{\red}$ is isomorphic over $K$ to a variety in $\Pb^{N-j}$ underlying a closed subscheme defined by equations of type $\bm^j$, and is non-empty whenever $N \geq |\bm^j|+j$.
\end{theorem}

We will also consider the following incidence variety:

\begin{definition}\label{def:total_polar_space}
The \emph{total $j$-polar space} of a variety $X$ is the space of all $j$-planes in $X$ equipped with an ordered spanning set: 
\[ \mcc^j(X) \coloneqq \left\{ (x_0,\dots,x_j,\Lambda) \in X^{j+1} \times \mcf(j,X) \mid \Lambda = \Lambda(x_0,\dots,x_j) \right\}. \]
Forgetting particular elements in these spanning tuples defines families of maps between total polar spaces. We write:
\[
\begin{aligned}
    p_j: \mcc^j(X) & \to \mcc^{j-1}(X) \\
    (x_0,\dots,x_j,\Lambda) & \mapsto (x_0,\dots,x_{j-1},\Lambda(x_0,\dots,x_{j-1})).
\end{aligned}
\]
\end{definition}

\begin{proposition}[{\cite[Proposition~3.3]{GGW2025}}]\label{prop:fundamental_span}
Let $X \subseteq \Pb^N$ be an intersection of hypersurfaces and take $j \in \Nb_{\geq 1}$. Consider the incidence correspondence
\[
\begin{tikzcd}
    & \mcc^{j}(X) \arrow[dl, "p_{j}"'] \arrow[dr, "q_{j}"] & \\
    \mcc^{j-1}(X) & & \mcf(j,X)
\end{tikzcd}
\]
where $p_j$ is as in Definition \ref{def:total_polar_space} and $q_{j}(x_0,\dots,x_{j},\Lambda) = \Lambda$. Then: 
\begin{enumerate}[label=(\alph*)]
\item The map $q_{j}$ is a Zariski locally trivial fiber bundle with fiber
\[ q_{j}^{-1}(\Lambda) = \PGL(\Lambda)/T(\Lambda), \]
where $T(\Lambda)$ is a torus conjugate to a subgroup of diagonal matrices, viewed as the stabilizer of an ordered spanning tuple $(x_0,\dots,x_{j})$ defined up to conjugacy in $\PGL(\Lambda)$.
\item The fiber $p_{j}^{-1}(x_0,\dots,x_{j-1},\Lambda)$ is given by the complement of the trivial sub-bundle $\Lambda$ in the restriction of the tautological bundle over $\grass(j,N)$ to $\mcf_\Lambda(j,X)$. 
\end{enumerate}
\end{proposition}

In particular, we will use:

\begin{lemma}\label{lemma:density_from_F_to_C}
Let $K$ be a $k$ field, let $\mce$ be a saturated class of accessory irrationalities, let $X$ be an intersection of hypersurfaces, and take $j \in \Nb$. If $\mcf(j,X)(K^{\mce}) \subset \mcf(j,X)$ is dense, then $\mcc^j(X)(K^{\mce}) \subset \mcc^j(X)$ is dense.
\end{lemma}
\begin{proof}
Since rational points of the fibers of the surjection $q_j: \mcc^j(X) \to \mcf(j,X)$ are dense over the field of definition of $x_0,\dots,x_j$, we see that $\mcc^j(X)(K^\mce) \subset \mcc^j(X)$ is dense whenever $\mcf(j,X)(K^{\mce}) \subset \mcf(j,X)$ is.
\end{proof}

\section{Special points via total \texorpdfstring{$j$}{j}-polar varieties}\label{sec:special_points_on_variety_of_j-planes}

\subsection{Obliteration lemmas}

We will now describe the obliteration algorithm of Sylvester \cite{Sylvester1887}, following the geometric approach in \cite{HS2023}.

\begin{definition}
    For $\mce$ a class of accessory irrationalities, we define 
    \[ f^{\mce}: \Nb \times \Mt \to \Nb \cup \{\infty\} \]
    by
    \begin{align*}
        f^{\mce}_j(\bm) \coloneqq \min\{ N \in \Nb \mid & \text{ for any } K/k, \text{ any } X \subset \Pb^N_K \text{ of type } \bm, \text{ and any } i \le j,  \\
        & \quad \mcf(j,X)(K^{\mce}) \subseteq \mcf(j,X) \text{ is nonempty and dense, and }\\
        & \quad p_{i+1} \circ \cdots \circ p_j\colon \mcc^j(X)\to \mcc^i(X) \text{ is dominant}\},
    \end{align*}
or $\infty$ if no such $N$ exists. 
\end{definition}
 
Since all our sequences have finite support, we will often write
\[ f^\mce_j(m_1,\dots,m_n) \coloneqq f^\mce_j(\bm), \] 
where $m_n > 0$ and $m_i = 0$ for all $i>n$. Moreover, when $\mce = \RD_{k}^{\leq \ell}$, we simply write $f^{\ell} \coloneqq f^\mce$.
We note that passing to hyperplanes (as in Remark \ref{rem:passing_to_hyperplanes}) allows us to omit linear terms from our calculations:
\begin{equation}\label{eqn:linear-forms-are-easy}
    f_j^{\mce}(m_1,m_2,\dots,m_n) = f_j^{\mce}(0,m_2,\dots,m_n) + m_1.
\end{equation}

\begin{remark}
    The class $\mce = \Sol$ is studied in \cite{GGW2025}, where it is simply denoted $f_j(m_2,m_3,m_4)$ (note that their indexing begins in degree $2$). We will establish a more general family of bounding lemmas here, where we are able to use the obliteration lemmas to greater effect when replacing $\Sol$ with the classes $\RD_k^{\leq \ell}$.
\end{remark}

In general, we are unable to compute exact values for $f^{\mce}_j$ except for some easy cases, e.g., 
\[ f^{\mce}_j(m_1) = j \quad \text{ and } \quad f^{\mce}_0(m_1) = m_1, \]
which are the minimal dimensions for $\grass(j,\Pb^N)$ and the intersection of $m_1$ generic linear forms, respectively, to be non-empty. The goal of this section is, given $\bm$, to develop the machinery needed in order to give upper bounds on $f^\ell_j(\bm)$.

\begin{remark}
The methods developed here only produce bounds on $f^\ell_j(\bm)$ when $\RD_k(n) \leq \ell$. Given a strict lower bound $\ell < \RD_k(n)$, alternate techniques (e.g., variations on Brauer diagonalization) can still be used to produce level $\ell$ points on such intersections, though we are unaware of an approach that will guarantee \emph{density} of such points.
\end{remark}

Firstly, we establish various forms of monotonicity enjoyed by the family of functions $f^{\mce}_j$. In what follows, we write $\leq$ for the entrywise partial order on $\Mt$.

\begin{proposition}\label{prop:monotonicity}
Let $\mce$ be a class of accessory irrationalities, let $\bm \in \Mt$, and let $j \in \Nb$. Then:
\begin{enumerate}[label=(\alph*)]
    \item $|\bm| \leq f^{\mce}_0(\bm)$
    \item $f^{\mce}_j(\bm)\le f^{\mce}_{j+1}(\bm)$.
    \item For any $\ba \in \Mt$ with $\ba \leq \bm$, we have $f^{\mce}_j(\ba)\le f^{\mce}_j(\bm)$.
    \item For any class $\mce'$ of accessory irrationalities with $\mce'\subseteq\mce$, we have $f^{\mce}_j(\bm)\leq f^{\mce'}_j(\bm)$.
\end{enumerate}
\end{proposition}
\begin{proof}
Part (a) is clear because, for $N < |\bm|$, the generic variety of type $\bm$ in $\Pb^N$ is empty. 

For (b), suppose that $f^{\mce}_{j+1}(\bm) \leq N$. Let $K/k$ be a field and let $X \subseteq \Pb^N_K$ be an intersection of type $\bm$. We will make use of the maps
\begin{equation}\label{eqn:inductive_span} 
\begin{tikzcd}
    \mcc^{j+1}(X) \arrow[r, "p_{j+1}"] \arrow[d, "q_{j+1}"'] & \mcc^j(X) \arrow[d, "q_j"] \\
    \mcf(j+1,X) & \mcf(j,X)
\end{tikzcd}
\end{equation}
described in Proposition \ref{prop:fundamental_span}. We have that $\mcf(j+1,X)(K^{\mce}) \subset \mcf(j+1,X)$ is dense and, by Lemma \ref{lemma:density_from_F_to_C}, so is $\mcc^{j+1}(X)(K^\mce) \subset \mcc^{j+1}(X)$. Since $p_{j+1}$ is dominant by the definition of $f^{\mce}$ and $q_j: \mcc^{j}(X) \to \mcf(j,X)$ is always surjective, we conclude that $\mcf(j,X)(K^{\mce}) \subset \mcf(j,X)$ is dense.

For (c), suppose $\ba \leq \bm$. By induction, it suffices to consider the case where $\bm$ is obtained by adding a single form of degree $d$ to $\ba$, which we write as $\bm = \ba + \be_d$. We will show that
\[
    f^{\mce}_j(\ba) < f^{\mce}_j(\ba + \be_1) \leq f^{\mce}_j(\ba+\be_d),
\]
assuming these values are finite.
The first inequality follows immediately from the definitions, e.g., \eqref{eqn:linear-forms-are-easy}. For the second inequality, let $K/k$ be a field extension, let $X \subseteq \Pb^N_K$ be an arbitrary intersection of type $\ba+\be_1$, and let $g$ denote one of its linear forms.
We can associate to $X$ an intersection $X'$ of type $\ba+\be_d$ by replacing $g$ with $g^d$. Since a form $g^d$ vanishes identically on a linear subspace $\Lambda$ if and only if $g$ vanishes on $\Lambda$, the $j$-planes contained in $X'$ are exactly those contained in $X$. Thus, $\mcf(j, X)(L) = \mcf(j, X')(L)$ for any field extension $L/K$. Since $X'$ is of type $\ba+\be_d$, the dimension $N$ required to guarantee $\mcf(j, X')(K^\mce) \subset \mcf(j,X')$ is dense is bounded by $f^{\mce}_j(\ba+\be_d)$, which thus also suffices for $\mcf(j,X)$.

Lastly, part (d) follows because $\mce'\subseteq\mce$ means that $K^{\mce'} \into K^{\mce}$.
\end{proof}

We now prove our two main results towards bounding $f_0^\ell(\bm)$.

\begin{lemma}[Linear subspaces via polar cones]\label{lemma:planes_from_pts} Fix $j \geq 0$, $\bm \in \Mt$, and a base field $k$ with $\characteristic(k) = 0$ or $\characteristic(k) > n$. If $\mce$ is a saturated class of accessory irrationalities closed under extensions, then
\begin{equation}\label{eqn:inductive_hypothesis} 
f^{\mce}_j(\bm) \leq f^{\mce}_0\left(\bm^j\right) + j.
\end{equation}
\end{lemma}

\begin{proof}
The argument follows \cite[Lemma 4.5]{GGW2025}, replacing $\Sol$ with $\mce$. We induct on $j$, where the base case $j=0$ is tautological. For the inductive step, we suppose that $K/k$ is a field extension and $X \subseteq \Pb^N_K$ is an intersection of hypersurfaces of type $\bm$. Assume the inductive hypothesis \eqref{eqn:inductive_hypothesis} together with the premise
\begin{equation}\label{eqn:inductive_premise}
f^{\mce}_0\left(\bm^{j+1}\right) + j+1 \leq N, 
\end{equation}
and we will deduce that $\mcf(j+1,X)(K^{\mce}) \subset \mcf(j+1,X)$ is dense.

We again make use of the maps in Equation \eqref{eqn:inductive_span}. Firstly, because $\bm^j \leq \bm^{j+1}$ is clear, by Proposition \ref{prop:monotonicity}(c) the inductive hypothesis applies and so $\mcf(j,X)(K^{\mce}) \subset \mcf(j,X)$ is dense. By Lemma \ref{lemma:density_from_F_to_C}, we also know that $\mcc^j(X)(K^\mce) \subset \mcc^j(X)$ is dense. 
Next, we consider the fibers of $p_{j+1}: \mcc^{j+1}(X) \to \mcc^{j}(X)$; let $\Lambda \in \mcf_{\Lambda}(j,X)$. 
By Remark \ref{rem:passing_to_hyperplanes} and Theorem \ref{thm:non-empty_fano_dim}(a), we regard $\mcf_{\Lambda}(j+1,X)^{\red}$ as type $\bm^{j+1}$ in projective $(N - (j+1))$-space. 
We know that $|\bm^{j+1}| \leq f^{\mce}_0(\bm^{j+1})$ by Proposition \ref{prop:monotonicity}(a), so with \eqref{eqn:inductive_premise} we conclude 
\[
    \left|\bm^{j+1}\right| + j+1 \leq N.
\]
Theorem \ref{thm:non-empty_fano_dim}(b) applies, so we know $\mcf_{\Lambda}(j+1,X) \not = \emptyset$ and thus $p_{j+1}$ is surjective. 
Indeed, since the type of $\mcf_{\Lambda}(j+1,X)^{\red}$ in $\Pb^{N-j-1}_K$ is $\bm^{j+1}$ as in Theorem \ref{thm:non-empty_fano_dim}(a), our assumption \eqref{eqn:inductive_premise} means that
\[
    \mcf_{\Lambda}(j+1,X)^{\red}(K^{\mce}) \subset \mcf_{\Lambda}(j+1,X)^{\red}
\] 
is dense for all $(x_0,\dots,x_j,\Lambda) \in \mcc^j(X)(K^{\mce})$. By appealing to the density of $K^{\mce}$-points in $\PGL(\Lambda)/T(\Lambda)$, the fibers of $p_j$ have dense $K^{\mce}$-points and thus $\mcc^{j+1}(X)(K^\mce) \subset \mcc^{j+1}(X)$ is dense. Finally, by dominance of $q_{j+1}$, we see that $\mcf(j+1,X)(K^{\mce}) \subset \mcf(j+1,X)$ is dense.
\end{proof}

\begin{lemma}[Obliteration via linear subspaces]\label{lemma:pts_from_planes} Let $\bm \in \Mt$ with $\bm = (m_1,\dots,m_n)$ and $\RD_k(n) \leq \ell$. 
Suppose further that $\ba = (a_1,\dots,a_n) \leq \bm$ entrywise and $\RD_k\big( 2^{a_2} \cdots n^{a_n} \big) \leq \ell.$ Then we have
\[ f^\ell_0(\bm) \leq f^\ell_{|\ba|}(\bm-\ba). \]
\end{lemma}
\begin{proof}
We adapt the argument of \cite[Lemma 4.6]{GGW2025} to the class $\RD_k^{\leq \ell}$ of accessory irrationalities. We employ the incidence correspondence
\[ \begin{tikzcd}
& \mci\coloneqq\{ (y,\Lambda) \in Y \times \mcf(j,Z) \, | \, y \in \Lambda \} \arrow[swap]{dl}{p} \arrow{dr}{q} \\
Y\cap Z && \mcf(j,Z),
\end{tikzcd} \]
where $K/k$ is a field, $Y, Z \subseteq \Pb^N_K$ are projective varieties, $p$ forgets the $j$-plane, and $q$ forgets the point. Note that $p^{-1}(y) \cong \mcf_y(j,Z)$, while $q^{-1}(\Lambda) = Y \cap \Lambda$. Provided that $\RD_k( \deg Y ) \leq \ell$, the fiber over $\Lambda \in \mcf(j,Z)(K^{(\ell)})$ will have dense $K^{(\ell)}$ points. If we also have that $\mcf(j,Z)(K^{(\ell)}) \subset \mcf(j,Z)$ is dense, then $\mci(K^{(\ell)}) \subset \mci$ must be dense; further, if $p$ is surjective (in particular, dominant), it follows that $(Y \cap Z)(K^{(\ell)}) \subset Y \cap Z$ will also be dense.

We will show that, if $f^\ell_{|\ba|}(\bm-\ba) \leq N$, then any intersection $X \subseteq \Pb^N_K$ of type $\bm$ has dense $K^{(\ell)}$-points. Let $Y$ be the intersection of the hypersurfaces corresponding to the forms in $\ba$, set $Z$ to be the intersection of the remaining forms so that $X = Y \cap Z$, and write $j = |\ba|$. By the assumption $f^\ell_{j}(\bm-\ba) \leq N$, we know that $\mcf(j,Z)$ has dense $K^{(\ell)}$-points. We also conclude the surjectivity of $p$ (by the definition of $f^\ell_j$). so it follows that $K^{(\ell)}$-points are dense in $X = Y \cap Z$.
\end{proof}

\begin{remark}\label{remark:lower_bound}
As per Example \ref{example:e-closures} and Proposition \ref{prop:monotonicity}(d), $f^\ell_j(\bm)$ is monotonically decreasing in $\ell$. On the other hand, $f^\ell_0(\bm)$ is bounded below by $|\bm|$. Indeed, if $\RD_k(2^{m_2} \cdots n^{m_n}) \leq \ell$ then $f^\ell_0(\bm) = |\bm|$, i.e. the minimum dimension to guarantee a non-empty intersection. For $\ell$ sufficiently large to resolve the inequalities via a single application of Lemmas \ref{lemma:pts_from_planes} and \ref{lemma:planes_from_pts},
\[ f^\ell_j(\bm) = f^{\ell}_0\left(\bm^j\right) + j = f^{\ell}_{|\bm^j|}(0) + j = |\bm^j| + j \]
recovers the type of $\mcf_\Lambda(j,X)^{\red}$ described in Theorem \ref{thm:non-empty_fano_dim}.
\end{remark}

These results allow us to obtain upper bounds on $f^\ell_j(\bm)$ in terms of other expressions of $f^\ell$. We will henceforth refer to a single application of Lemmas \ref{lemma:pts_from_planes} and \ref{lemma:planes_from_pts} as ``obliterating $\bm$ by $\ba$'':

\begin{corollary}[Obliterating $\bm$ by $\ba$]\label{cor:obliteration}
Let $\bm = (m_1,\dots,m_n) \in \Mt$ with $\RD_k(n) \leq \ell$. Whenever $\ba = (a_1,\dots,a_n) \in \Mt$ has $\ba \leq \bm$ and $\RD_k\big( 2^{a_2} \cdots n^{a_n} \big) \leq \ell$, we have
\[
    f^\ell_0(\bm) \leq f^\ell_0\left( (\bm-\ba)^{|\ba|} \right) + |\ba|,
\]    
\end{corollary}

\subsection{Bounds on \texorpdfstring{$f^\ell_j(\bm)$}{flj(m)}}\label{subsec:bounds_on_f}

We remind the reader that, to date, no non-trivial lower bounds on $\RD_k(n)$ have been exhibited. To speak concretely about the process of improving existing estimates, we make the following explicit:

\begin{definition}\label{def:bound_on_RD}
Let $k$ be a field. A function $R: \Nb_{\geq 1} \to \Nb_{\geq 1}$ is called a \emph{bound on resolvent degree} if
\begin{enumerate}[label=(\alph*)]
    \item $\RD_k(n) \leq R(n) \leq n$ for all $n \in \Nb_{\geq 1}$,
    \item $R$ is weakly increasing,
    \item\label{def:hamilton_condition} For all $n \in \Nb_{\geq 1}$, $R(n+1) \leq R(n) + 1$.
\end{enumerate}
\end{definition}

We refer to the last of these as the \emph{Hamilton condition} (see Lemma \ref{lemma:hamilton_properties} and Appendix \ref{appendix:elementary_speed-ups} on the use of this condition). We note that all upper bounds on $\RD_k(n)$ of which we are aware in the literature satisfy this condition. However, we are not able to resolve the following:

\begin{problem}
    Is $\RD_k(n+1)\le \RD_k(n)+1$ for all $n$?
\end{problem}

Given how rapidly $\bm \mapsto \bm^j$ grows, a natural question is whether this process terminates---certainly, poor obliteration strategies (in particular, obliterating lower order terms while neglecting higher order ones) will produce arbitrarily large types in the bounding terms if continued indefinitely. In order to systematize our choices of obliteration, we introduce the following:

\begin{definition}
Given $\bm = (m_1,\dots,m_n)$ and a bound on resolvent degree $R$ with $R(n) \leq \ell$, a \emph{valid obliteration} (of $\bm$) is any $\ba \in \Mt$ with $\ba \leq \bm$ satisfying 
\[
    R\big( 2^{a_2} \cdots n^{a_n} \big) \leq \ell.
\]
With respect to a given bound $R$ and level $\ell$, the set of valid obliterations of $\bm$ is exactly those $\ba$ for which Corollary \ref{cor:obliteration} can be applied.
\end{definition}

For the remainder of the paper, we refer to the process of ``bounding $f^\ell_j(\bm)$'' to mean applying Corollary \ref{cor:obliteration} until a base case $f^\ell_j(m_1)=m_1+j$ is reached. Some calculations require so many intermediary steps that we do not record them, referring the reader to Appendix \ref{appendix:computing_f} and \cite{GomezGonzales2025} for details on reproducing the results. Practically, this process is carried out using a particular bound on resolvent degree. Moreover, at each step, we employ a top-down greedy obliteration strategy: \emph{We choose the valid obliteration $\ba \leq \bm$ subject to the condition that $a_i = 0$ whenever there exists some $j > i$ with $a_j < m_j$ and which sequentially maximizes its entries from right to left}. In other words, we never obliterate a lower-order term unless all higher-order terms have been completely eliminated. We justify this approach before proceeding.

\begin{example}\label{example:too_greedy_obliteration}
    Consider the intersection of a quadric $Q$, a cubic $C$, and a quartic $T$ in the case $\ell=4$. In order to make apparent what is being obliterated, we explicitly apply Lemmas \ref{lemma:pts_from_planes} and \ref{lemma:planes_from_pts} separately in the calculations to follow.

    First we demonstrate the method outlined above:
    \[
    \begin{aligned}
        f_0^4(0,1,1,1) & \leq f_1^4(0,1,1) \leq f_0^4(2,2,1) + 1 \leq f_2^4(2,1) + 1 \\
        & \leq f_0^4(4,1) + 3 \leq f_1^4(4)+3 \leq f_0^4(4) + 4 = 8,
    \end{aligned}
    \]
    where the first, third, and fifth inequalities use that $\RD_k(4), \RD_k(6), \RD_k(2) \leq 4$, respectively. We also implicitly rely on the fact that, while we do not have lower bounds for $\RD_k(12)$, our best current upper bound (for $k=\Cb$) is $\RD_k(12) \leq 7$, and hence we cannot eliminate any additional terms in the first two obliterations. 
    
    We also demonstrate a greedier approach which obliterates all forms available:
    \[
    \begin{aligned}
        f_0^4(0,1,1,1) & \leq f_2^4(0,0,1) \leq f_0^4(3,2,1) + 2 \leq f_2^4(3,1) + 2 \\
        & \leq f_0^4(5,1) + 4 \leq f_1^4(5)+4 \leq f_0^4(5) + 5 = 10.
    \end{aligned}
    \]
    In particular, at the first obliteration, the greedy approach obliterates both $Q$ and $T$, since $\RD_k(8) \leq 4$, at the cost of passing to planes rather than lines, but this leads to a worse bound.
\end{example}

Indeed, the approach outlined here is an effective greedy heuristic, though not always optimal. We demonstrate this in two parts, together with an example in which a different strategy improves the bound. First, we establish that obliterating by batches produces better bounds on $f^\ell_0(\bm)$ when compared with sequential obliterations that eliminate the same forms.

\begin{proposition}\label{prop:batching_obliterations}
Let $\bm = (m_1,\dots,m_n), \ba, \bb \in \Mt$. If $\ba + \bb \leq \bm$, $\bb \leq (\bm - \ba)^{|\ba|}$, and
\[ \RD_k(2^{a_2+b_2} \cdots n^{a_n+b_n}) \leq \ell, \]
then the obliteration algorithm from obliterating by $\ba+\bb$ yields a bound on $f_0^{\ell}(\bm)$ which is less than or equal to the bound from obliterating by $\ba$ followed by obliterating by $\bb$. 
\end{proposition}

\begin{proof}
We compute using Corollary \ref{cor:obliteration}. Obliterating in a single batch gives
\[
    f^\ell_0(\bm) \leq f^\ell_0\left( (\bm-(\ba + \bb))^{|\ba|+|\bb|} \right) + |\ba|+|\bb|,
\]
while obliterating sequentially yields the bound
\[
    f^\ell_0(\bm) \leq f^\ell_0\left( (\bm-\ba)^{|\ba|} \right) + |\ba| \leq f^\ell_0\left( \left( (\bm - \ba)^{|\ba|} - \bb \right)^{|\bb|} \right) + |\ba| +|\bb|.
\]
We know by Lemma \ref{lemma:batching_obliterations} that
\[
    (\bm - (\ba + \bb))^{|\ba|+|\bb|} \leq \left( (\bm - \ba)^{|\ba|} - \bb \right)^{|\bb|}.
\]
and so we conclude by monotonicity (Proposition \ref{prop:monotonicity}).
\end{proof}

Next, we generalize Example \ref{example:too_greedy_obliteration}, showing that obliterating the highest terms in $\bm$ produces better bounds.

\begin{proposition}\label{prop:skipping_rightmost}
Let $\bm = (m_1,\dots,m_n), \ba = (a_1, \dots, a_{n-1}) \in \Mt$. If $0 \not = \ba \leq \bm$, $0 \leq j < m_n$, and 
\[ \RD_k(2^{a_2} \cdots (n-1)^{a_{n-1}} n^j) \leq \ell, \]
then obliterating $j$ forms in degree $n$ using the obliteration algorithm produces a bound on $f_0^\ell(\bm)$ strictly less than that obtained by obliterating $j$ forms in degree $n$ together with the lower-degree terms specified by $\ba$.
\end{proposition}
\begin{proof}
We bound $f_0^\ell(\bm)$ in two ways, obliterating by $j \be_n$ versus $j \be_n + \ba$, which gives
\[
    f_0^\ell(\bm) \leq f_0^\ell\left((\bm-j \be_n)^j\right) + j \quad \text{ and } \quad
    f_0^\ell(\bm) \leq f_0^\ell\left((\bm-j \be_n-\ba)^{j+|\ba|}\right) + (j+|\ba|).
\]
By Proposition \ref{prop:skipping_rightmost_combinatorics}, we know that
\[
    (\bm - j \be_n )^j \leq (\bm - j \be_n - \ba)^{j+|\ba|}
\]
holds entrywise, and so we conclude by monotonicity and the fact that $j < j+|\ba|$.
\end{proof}

\noindent In particular, the $j=0$ case of the above demonstrates that obliterating lower order terms without first dealing with the highest order terms yields worse bounds. 

Taken together, Propositions \ref{prop:batching_obliterations} and \ref{prop:skipping_rightmost} justify our top-down greedy obliteration heuristic. Specifically, they establish that for any individual application of Corollary \ref{cor:obliteration} to a tuple $\bm$, choosing 
\begin{equation}\label{eqn:obliteration_strategy}
\begin{aligned}
    \ba \coloneqq \max\big\{ (a_1,\dots,a_n) \mid\ & (a_1,\dots,a_n) \text{ is a valid obliteration of } \bm \text{ and } \\
    & \quad \text{for any } i < d \leq n,\ a_i = 0 \text{ whenever } a_d < m_d \big\}, \\
\end{aligned}
\end{equation}
minimizes the resulting upper bound expression in that step. However, while this greedy selection is locally optimal, it does not guarantee global optimality over a sequence of obliterations. Towards this end, we give three examples (over $k = \Cb$) to illustrate the process as types grow more complicated before proving Theorem \ref{thm:polynomial_bounds_body}.

\begin{example}\label{example:non_greedy_better}
Consider three quadrics and one cubic at $\ell = 4$. Our strategy selects $\ba = (0,1,1)$ since $\RD(2 \cdot 3) \leq 2 \leq 4$, while the current best bound on $\RD(2^2 \cdot 3)$ is $7$. This yields:
\[
    f_0^4(0,3,1) \le f_2^4(0,2) \le f_0^4(4,2) + 2 \le f_2^4(4) + 2 = 4 + 2 + 2 = 8,
\]
where we use $\RD(2^2) = 1$ in the final obliteration. On the other hand, if we only eliminate the cubic via $\ba = (0,0,1)$, we have
\[
    f_0^4(0,3,1) \le f_1^4(0,3) \le f_0^4(3,3) + 1 \le f_3^4(3,0) + 1 = 3 + 3 + 1 = 7,
\]
using the fact that $\RD(2^3) \leq 4$. In this case, deferring the quadrics yields a better bound than our heuristic strategy.
\end{example}

Despite such examples, this greedy heuristic provides an effective, deterministic baseline for computing bounds on $f_0^\ell(\bm)$ as $\ell$ varies, avoiding the combinatorial explosion of broader state-space searches as types grow, which is needed for the proof of Theorem \ref{thm:new_bounds_on_H}. Finding the globally minimal bound across all possible sequences of obliterations remains an open question.

\begin{problem}\label{prob:optimal_obliteration}
Given $\bm \in \Mt$, a level $\ell$, and bound on resolvent degree $R$, determine an efficient algorithm to compute the globally optimal sequence of valid obliterations minimizing the upper bound on $f_0^\ell(\bm)$ under Corollary \ref{cor:obliteration}. In particular, under what conditions does the top-down greedy obliteration \eqref{eqn:obliteration_strategy} strategy coincide with the global optimum?
\end{problem}

\begin{example}\label{example:computations_for_nonic} 
    We compute a bound on $f^\ell_0(0,1,1,1)$ for all $\ell$ using our greedy strategy:

    \[
    \begin{tikzcd}[row sep=normal, column sep=normal]
        & & f^\ell_0(0,1,1,1) \arrow[ddl, "\ell < \RD(12)"'] \arrow[dd, "\RD(12) \leq \ell < \RD(24)" description] \arrow[ddr, "\ell \geq \RD(24)"] & \\
        & & & \\
        & \leq f^\ell_1(0,1,1) \arrow[d] & \leq f^\ell_2(0,1) \arrow[d] & = f^\ell_3(0) = 3 \\
        & \leq f^\ell_0(2,2,1)+1 \arrow[dl, "\ell < \RD(6)"'] \arrow[d, "\RD(6) \leq \ell"] & \leq f^\ell_0(2,1)+2 = 5 & \\
        \leq f^\ell_1(2,2)+1 \arrow[d] & \leq f^\ell_2(2,1)+1 \arrow[d] & & \\
        \leq f^\ell_0(4,2)+2 \arrow[d] & \leq f^\ell_0(4,1)+3 = 8 & & \\
        \leq f^\ell_6(0)+2 = 8 & & & \\
    \end{tikzcd}
    \]
    
    The inequalities above reflect an implicit algorithm for obtaining level $\ell$ points on intersections of hypersurfaces. Consider the case\footnote{We emphasize that such an $\ell$ might not exist---there are no known non-trivial lower bounds on $\RD(n)$.} where $\RD(12) \leq \ell < \RD(24)$. Suppose we have a quadric $Q$, a cubic $C$, and a quartic $T$ in $\Pb^N$, and we want to find a level $\ell$ point on $Q \cap C \cap T$. To do so, we will explicitly construct a $2$-plane $\Lambda \subset Q$, which we then intersect with $C \cap T$. Because $\Lambda \cong \Pb^2$, this amounts to finding a point in the intersection of two plane curves (of degrees 3 and 4), which can be done at the cost of solving a degree 12 polynomial (e.g., via the resultant).

    In order to construct $\Lambda \in \mcf(2, Q)$, we build it up dimension by dimension: we first find a point $y_1 \in Q$ (by solving a quadratic), and then a line $L \subset Q$ containing $y_1$ (by solving another quadratic, together with some additional linear constraints). We then consider $\mcf_{L}(2,Q)^{\red}$, which is isomorphic to a variety of type $(2,1)$ in $\Pb^{N-2}$. We proceed by choosing a $3$-plane $\Lambda' \in \mcf(3,\Pb^{N-2})$, which is always possible for $N-2 \geq 3$, and intersecting it with $\mcf_{L}(2,Q)^{\red}$, ultimately finding our plane $\Lambda$ by solving one final quadratic.

    From this perspective, the greedy obliteration in Example \ref{example:too_greedy_obliteration} leads to a worse bound because the type of $\mcf_{L}(2,C)$ exceeds the type of $\mcf_{y}(1,Q\cap C)$.
\end{example}

With our strategy \eqref{eqn:obliteration_strategy} for obliteration established and unambiguous, we will henceforth apply Corollary \ref{cor:obliteration} in a single step (i.e., without writing the inequalities of Lemmas \ref{lemma:planes_from_pts} and \ref{lemma:pts_from_planes} separately).

\begin{example}\label{example:many_accessory_irrationalities}
Computing a bound on $f^{70959641905151979}_0(\overbrace{1,\dots,1}^{21})$ by the greedy heuristic---which arises in the course of proving Theorem \ref{thm:new_bounds_on_H} in Section \ref{subsec:new_bounds}---involves $1029976$ applications of Corollary \ref{cor:obliteration} to arrive at the base case $f^\ell_0(m_1) = m_1$. We summarize the process here by recording the total number of applications of the Corollary over the course of evaluating the bound, along with the associated length $n$ of the tuple at each stage.
\begin{center}
    \begin{tabular}{ccc}
        \toprule
        $n$ & Corollary \ref{cor:obliteration} \\
        \midrule
        $21$ & 1 \\
        $6$ & 1 \\
        $4$ & 5 \\
        $3$ & 302 \\
        $2$ & 1029667 \\
        \bottomrule
    \end{tabular}
\end{center}
\end{example}

Given the growth of $\bm^j$ and the number of lower-order terms which  accumulate over the course of evaluating a bound on $f^\ell_0(\bm)$, it is difficult to know from the start which obliterations to expect over the course of a calculation (or how many). In Proposition \ref{prop:l-points_on_quadrics} of Appendix \ref{appendix:elementary_speed-ups} we provide the closed-form expression for the bound on $f^\ell_0(m_1,m_2)$ in terms of how $R(2^b)$ grows, but in general we do not have an effective method for writing down sharp (with respect to the methods of this paper) bounds on $f^\ell_0(\bm)$ outside of explicit iterative computation with specific $\bm$, $\ell$, and $R$.

With this in mind, we can establish a \emph{coarse} upper bound on $f^\ell_j(\bm)$ by obliterating using only lines, i.e., by only obliterating a single form at a time, from which we deduce Theorem \ref{thm:polynomial_bounds_body}. See Appendix \ref{appendix:computing_f}, especially Example \ref{example:coarse_bound_is_coarse}, for a comparison between this upper bound and those derived by the greedy heuristic \eqref{eqn:obliteration_strategy}.

\begin{lemma}[Reduction of degree by lines]\label{lemma:coarse_upper_bound}
Let $\bm = (m_1, \dots, m_n) \in \Mt$, let $k$ be a base field with $\characteristic(k) = 0$ or $\characteristic(k) > n$, and let $\ell \geq \RD_k(n)$. Then we have
\[
    f^\ell_0(\bm) \leq f^\ell_0\bigg( \bm^{m_n} - \sum_{i=1}^{m_n} \be_n^i \bigg) + m_n.
\]
The argument of $f^\ell_0$ on the right-hand side is supported up to degree $n-1$, arising from $m_n$ successive obliterations by a single form in degree $n$. In degree $d < n$, the argument is given by
\[
    \sum_{i=d}^{n-1} \binom{m_n+i-d-1}{i-d} m_i + (n-d) \binom{m_n+n-d-1}{n-d+1}.
\]
\end{lemma}

\begin{proof}
We show a more general result, given by sequentially obliterating $s$ forms in degree $n$:
\[
    f^\ell_0(\bm) \leq f^\ell_0\bigg( \bm^s - \sum_{j=1}^s \be_n^j \bigg) + s. 
\]
We proceed by induction on $s$, where the base case $s = 0$ is tautological. For the inductive step, applying Corollary  \ref{cor:obliteration} to the right-hand side of the above inequality, we see that
\[
\begin{aligned}
    f^\ell_0(\bm) \leq f^\ell_0\bigg( \bm^s - \sum_{j=1}^s \be_n^j \bigg) + s & \leq f^\ell_0\bigg( \big(\bm^s - \sum_{j=1}^s \be_n^j\big)^1 - \be_n^1 \bigg) + (s+1) \\
    & = f^\ell_0\bigg(\bm^{s+1} - \sum_{j=1}^{s+1} \be_n^j \bigg) + (s+1),
\end{aligned}
\]
which completes the first claim. The $d$th entry of $\bm^{m_n} - \sum_{j=1}^{m_n} \be_n^j$ for $1 \leq d < n$ is given by Lemma \ref{lemma:j-endomorphism}(b):
\[
    \sum_{i=d}^{n} \binom{m_{n}+i-d-1}{i-d} m_i - \sum_{j=1}^{m_n} \binom{j+n-d-1}{n-d}
\]
By Chu's theorem\footnote{This is also called the hockey-stick identity; see, e.g., Theorem~1.5.2 of \cite{Merris2003}.} we have
\[
    \sum_{i=d}^{n} \binom{m_n+i-d-1}{i-d} m_i - \binom{m_n+n-d}{n-d+1}.
\]
Taking the difference between the pure $m_n$ terms of these expansions gives
\[
\begin{aligned}
    \binom{m_n+n-d-1}{n-d} m_n - \binom{m_n+n-d}{n-d+1} & = \frac{(m_n+n-d-1)!}{(n-d+1)! (m_n-1)!} \big( (n-d+1)m_n - (m_n+n-d) \big) \\
    & = \frac{(m_n+n-d-1)!}{(n-d+1)! (m_n-1)!} (n-d)(m_n-1) \\
    & = (n-d) \binom{m_n+n-d-1}{n-d+1}.
\end{aligned}
\]
\end{proof}

Now, together with Remark \ref{remark:lower_bound}, we conclude Theorem \ref{thm:polynomial_bounds} in fuller detail:

\begin{theorem}\label{thm:polynomial_bounds_body}
Fix $n \geq 2$, a base field $k$ with $\characteristic(k) = 0$ or $\characteristic(k) > n$, and $\ell \geq \RD_k(n)$. Then there is a $p_n \in \Qb[m_1,\dots,m_n]$ of degree $2^{i-1}$ in $m_i$ and of total degree $2^{n-1}$ so that, for any $\bm = (m_1,\dots,m_n) \in \Mt$,
\[ |\bm| \leq f^\ell_0(\bm) \leq p_n(\bm). \]
Further, there is a polynomial $q_n \in \Qb[j,m_1,\dots,m_n]$ of degree $2^{i-1}$ in each $m_i$, degree $2^{n-2}$ in $j$, and total degree $2^{n-1}$ such that $p_n(\bm) = q_n(0,\bm)$ and
\[ 
    (j+1) + \sum_{i=1}^n \binom{j+i-1}{i-1} m_i \leq f^\ell_j(\bm) \leq q_n(j,\bm). 
\]
\end{theorem}

\begin{proof}
We have already established the desired lower bounds in Remark \ref{remark:lower_bound}. We proceed with the remainder of the proof by induction on $n$. For the base case $n=2$, we apply Lemma \ref{lemma:coarse_upper_bound} to $\bm = (m_1, m_2)$ to see
\[ 
    f_0^\ell(\bm) \leq f_0^\ell\left( m_1 + \binom{m_2}{2} \right) + m_2 \leq m_1 + \frac{1}{2}m_2(m_2+1) \eqqcolon p_2(\bm), 
\]
which has the desired degrees. Next, applying Lemma \ref{lemma:planes_from_pts} gives
\[ 
    f_j^\ell(\bm) \leq f_0^\ell(\bm^j) + j \leq p_2(\bm^j) + j = m_1 + \frac{1}{2} (m_2+1) (m_2+2j) \eqqcolon q_2(j,\bm), 
\]
since $\bm^j = (m_1+j m_2, m_2)$. This establishes the base case.

Suppose the claims hold for some $n \geq 2$ and let $\bm = (m_1,\dots,m_{n+1}) \in \Mt$. By Lemma \ref{lemma:coarse_upper_bound}, we have
\[
    f^\ell_0(\bm) \leq f^\ell_0\bigg(\bm^{m_{n+1}} - \sum_{i=1}^{m_{n+1}} \be_{n+1}^i\bigg) + m_{n+1},
\]
where, for any $1 \leq d \leq n$, the $d$th coordinate in the argument of $f^\ell_0$ on the right-hand side is
\[
    \sum_{i=d}^{n} \binom{m_{n+1}+i-d-1}{i-d} m_i + (n+1-d) \binom{m_{n+1}+n-d}{n-d+2}.
\]
In the $d$th entry, the variable $m_i$ for $i \leq n$ occurs with degree $1$; on the other hand, $m_{n+1}$ occurs with degree $n-d+2$. We define
\[
    p_{n+1}(\bm) \coloneqq p_n\bigg(\bm^{m_{n+1}} - \sum_{i=1}^{m_{n+1}} \be_{n+1}^i\bigg) + m_{n+1},
\]
so that $m_d$ appears in $p_{n+1}$ with degree at most the degree of the $d$th variable in $p_n$, which is $2^{d-1}$, and $m_{n+1}$ contributes to $p_{n+1}$ with degree $(n-d+2)2^{d-1}$ for each $d$. The latter degree is maximized when $d=n$ as $2 \cdot 2^{n-1} = 2^n$, so we conclude $p_{n+1}$ has the desired degree in $m_{n+1}$. For total degree, the highest total degree term in the $d$th entry is purely in $m_{n+1}$ with degree $n-d+2$, which again gives the total degree of $p_{n+1}$ as $2^n$. 
Lastly, we define $q_{n+1}$ by applying Lemma \ref{lemma:planes_from_pts} directly:
\[
    f^\ell_j(\bm) \leq f^\ell_0(\bm^j) + j \leq p_{n+1}(\bm^j) + j \eqqcolon q_{n+1}(j, \bm).
\]
The $d$th coordinate of $\bm^j$ contains $j$ with degree $n+1-d$, which contributes to $q_{n+1}$ with degree $(n+1-d)2^{d-1}$. This is maximized both when $d=n$ (giving $1 \cdot 2^{n-1}$) and $d=n-1$ (giving $2 \cdot 2^{n-2} = 2^{n-1}$), so the degree of $j$ in $q_{n+1}$ is $2^{n-1}$. Moreover, the total degree of the $d$th entry of $\bm^j$ is $n+2-d$, contributing to $q_{n+1}$ with total degree at most $2^n$. Thus $p_{n+1}$ and $q_{n+1}$ are polynomial bounds of the desired form.
\end{proof}

\section{Bounds on \texorpdfstring{$\RD(n)$}{RD(n)}}\label{sec:bounds_on_RD}

Given polynomials $f_1,\dots,f_s \in K[x_0,\dots,x_N]$, we write $\Vb(f_1,\dots,f_s)$ for the (scheme-theoretic) intersection of hypersurfaces $\bigcap_{j=1}^s \Vb(f_j)$ in $\Pb_K^N$. Further, for the remainder of the paper, we write $\RD(n) \coloneqq \RD_\Cb(n)$ and $K_n \coloneqq \Cb(a_1,\dots,a_n)$. 

\subsection{Known bounds}\label{sec:known_bounds}
Since Hamilton \cite{Hamilton1836} first proved that \[ \lim_{n \to \infty} n - \RD(n) = \infty, \] the literature for improving estimates on $\RD(n)$ has nearly always been phrased in terms of functions $H: \Nb_{\geq 1} \to \Nb_{\geq 1}$ such that $n - \RD(n) \geq r$ for all $n \geq H(r)$. As such, we formally introduce the following:

\begin{definition}\label{def:hamilton_function}
The \emph{Hamilton function} $\mch: \Nb_{\geq 1} \to \Nb_{\geq 1}$ is given by
\[ \mch(r) \coloneqq \min \{ n \in \Nb_{\geq 1}: n-\RD(n) \geq r \}. \]
\end{definition}

The classical formulas in radicals for quadratics, cubics, and quartics, together with Bring's \cite{Bring1786} reduction of the quintic, established that $\mch(1) = 2,$ $\mch(2) = 3,$ $\mch(3) = 4,$ and $\mch(4) = 5$. Later, Hilbert's sketch \cite{Hilbert1927}---which used the 27 lines on a smooth cubic surface to reduce the generic nonic equation---demonstrated that $\mch(5) \leq 9$. 

The significance of the Hamilton function is related to the central role of \emph{Tschirnhaus transformations}. These correspond to rational changes of variables designed to eliminate intermediate coefficients of a polynomial. Formally, they are isomorphisms over $K_n$ of the form
\[
    T: K_n[x]/(x^n+a_1 x^{n-1} + \cdots + a_n) \stackrel{\cong}{\longrightarrow} K_n[x]/(x^n+b_1 x^{n-1} + \cdots + b_n),
\]
where specific coefficients $b_i \in K_n$ are zero. The particular configuration of vanishing $b_i$ defines the type of the Tschirnhaus transformation. Here we focus on transformations where the first $r-1$ coefficients vanish: $b_1 = \cdots = b_{r-1} = 0$. Such an isomorphism reduces the generic degree $n$ polynomial to an algebraic function of fewer variables:
\begin{equation}\label{eqn:reduced_poly}
    (b_r,\dots,b_n) \mapsto \{ x \in \Cb: x^n + b_rx^{n-r} + \cdots + b_n = 0\}.
\end{equation} 
This can be further normalized by scaling to $x^n + b_r x^{n-r} + \cdots + b_{n-1} x + 1 = 0$, demonstrating that $\RD(n) \leq n-r$, and consequently, $\mch(r) \leq n$. 

As described in \cite{Coray1987, Dixmier1993}, any such change of variables can be expressed in terms of primitive elements as \[ \bar{x} \mapsto w_0 + w_1 \bar{x} + \cdots + w_{n-1} \bar{x}^{n-1}. \] Note that this assignment defines a Tschirnhaus transformation (i.e., an isomorphism) exactly when at least one of the $w_i \not = 0$ for $i > 0$. Moreover, because scaling the vector of coefficients $(w_0, \dots, w_{n-1})$ does not change the roots of the resulting polynomial, the space of all possible transformations is naturally parameterized by the projective space $\Pb_{K_n}^{n-1} \setminus \{ (1:0:\cdots:0) \}$. 

Each coefficient $b_i$ is a homogeneous polynomial of degree $i$ in the variables $w_0,\dots,w_{n-1}$ with coefficients in $K_n$. Therefore, finding a Tschirnhaus transformation that eliminates the first $r-1$ coefficients is equivalent to finding a point on the intersection of hypersurfaces
\[
    \tau_{1,\dots,r-1} \coloneqq \Vb\left(b_1,\dots,b_{r-1} \right) \subseteq \Pb_{K_n}^{n-1},
\]
where we mirror the notation of \cite{Wolfson2020,Sutherland2022} for consistency.

\begin{example}
The change of variables $\bar{x} \mapsto w_0 + w_1 \bar{x}$ sends $x^2 + a_1 x + a_2 = 0$ to
\[
    x^2 + (\underbrace{-2w_0 + a_1w_1}_{b_1})x + (\underbrace{w_0^2 - w_0 w_1 a_1 + w_1^2 a_2}_{b_2}) = 0.
\]
The point $(w_0:w_1) = (a_1:2) \in \Vb(b_1)$ corresponds to the usual substitution $x \mapsto x + \tfrac{a_1}{2}$ used to complete the square. 
\end{example}

In this geometric framework, guaranteeing the existence of a point on $\tau_{1,\dots,r-1} \setminus \{(1:0:\cdots:0)\}$ defined over some field extension $L/K_n$ of suitably low resolvent degree is sufficient to show $\mch(r) \leq n$.\footnote{For a fuller discussion, see Sections 3--4 of \cite{Wolfson2020}, especially Lemma 4.12.} Specifically, given such a point defined over $L \subseteq K_n^{(\ell)}$, we can eliminate $r$ variables from the generic polynomial at the cost of introducing an auxiliary tower of field extensions with resolvent degree at most $\ell$. This yields:
\[
    \RD(n) \leq \max(n-r,\ell).
\]

Recent improvements to Brauer's classical bounds for $\mch(r)$, such as those in \cite{Wolfson2020, Sutherland2021}, approach this problem by establishing the existence of an $(r-d-1)$-dimensional linear subspace $\Lambda \subseteq \tau_{1,\dots,d} \setminus \{ (1:0:\cdots:0) \}$ over an extension $L/K_n$ of controlled resolvent degree. Once this subspace is found, we intersect $\Lambda$ with the remaining hypersurfaces $\Vb(b_{d+1}, \dots, b_{r-1})$; this intersection is cut out by forms of degrees $d+1$ through $r-1$, meaning a valid Tschirnhaus transformation can be obtained at the cost of solving an auxiliary polynomial of degree $\tfrac{(r-1)!}{d!}$. This strategy gives bounds of the form:
\[
    \mch(r) \leq \tfrac{(r-1)!}{d!}+1.
\]

The methodologies used to establish the existence of the subspace $\Lambda$ depend heavily on the scale of $r$. For example, \cite{Sutherland2022} employs dimension-counting on moduli spaces to obtain asymptotic bounds. The methods developed in this paper appear to be most applicable to sharpening bounds in an ``intermediate range''---meaning values of $r$ that are larger than those addressed by classical geometric constructions (e.g., the $27$ lines on a cubic surface used for $r = 5$), but small enough that general asymptotic bounds remain suboptimal. In this regime, \cite{HS2023} relied on geometric constructions involving polar cones to identify such $(r-d-1)$-planes directly. Our generalized framework in this paper improves upon these constructions, applying the greedy obliteration strategy \eqref{eqn:obliteration_strategy} to yield sharper upper bounds for this intermediate range.

\subsection{New bounds}\label{subsec:new_bounds}

Recall (Definition \ref{def:bound_on_RD}) that a \emph{bound on resolvent degree} is a weakly increasing function $R: \Nb_{\geq 1} \to \Nb_{\geq 1}$ satisfying $\RD(n) \leq R(n) \leq n$ and also the \emph{Hamilton condition},
\[ R(n+1) \leq R(n)+1, \]
for all $n \in \Nb_{\geq 1}$.

\begin{definition}\label{def:hamilton_function_for_bound}
Let $R$ be a bound on resolvent degree. The \emph{Hamilton function for $R$} is the function $H_R: \Nb_{\geq 1} \to \Nb_{\geq 1} \cup \{\infty\}$ given by
\[ H_R(r) \coloneqq \min\{ n \in \Nb_{\geq 1}: n - R(n) \geq r \}, \]
or $H_R(r) \coloneqq \infty$ if $n - R(n) < r$ holds for all $n \in \Nb_{\geq 1}$.
\end{definition}

\begin{lemma}\label{lemma:hamilton_properties}
If $R$ is a bound on resolvent degree with Hamilton function $H_R$, then:
\begin{enumerate}[label=(\alph*)]
    \item $H_R$ is strictly increasing.
    \item If $H_R(r) \leq n < H_R(r+1)$, then $R(n) = n - r$.
    \item For any $r \in \Nb_{\geq 1}$ with $H_R(r) < \infty$, we have $R(H_R(r)) = H_R(r) - r$.
    \item For any $r \in \Nb_{\geq 1}$, $\mch(r) \leq H_R(r)$.
\end{enumerate}
\end{lemma}
\begin{proof}
Part (a) follows by definition. For part (b), we first observe that since $R$ is weakly increasing and satisfies the Hamilton condition, we always have $R(n) \leq R(n+1) \leq R(n)+1$. In particular, the difference $n-R(n)$ is also weakly increasing: applying these inequalities to $(n+1)-R(n+1)$ gives
\[ n - R(n) \leq (n+1) - R(n+1) \leq n+1 - R(n). \]
Part (b) follows from this, because the above combined with $H_R(r) \leq n < H_R(r+1)$ implies that $n - R(n) = r$. We get (c) by taking $n = H_R(r)$ and (d) by noting that the set being minimized over to compute $H_R(r)$ is contained in that used to compute $\mch(r)$.
\end{proof}

In particular, bounds on resolvent degree and their associated Hamilton function should be understood as carrying equivalent data; obtaining a tighter bound on $\mch$ is equivalent to improving bounds on $\RD$.

\begin{example}
The functions $F$, $G$, and $G'$ defined in \cite{Wolfson2020,Sutherland2021,HS2023}, respectively, are the Hamilton functions associated to bounds on resolvent degree established in each paper.
\end{example}

Given $R$ a bound on $\RD$, we can use it to compute bounds for $f^\ell_j(m_1,\dots,m_n)$ according to our obliteration strategy \eqref{eqn:obliteration_strategy} whenever $R(n) \leq \ell$. Our next result shows how to use such calculations to obtain a new bound on resolvent degree from $R$.

\begin{theorem}\label{thm:sharpening_bound_on_RD}
Fix $r \geq 4$ and let $R$ be a bound on $\RD$ with Hamilton function $H_R$. Then
\begin{equation}\label{eqn:sharpenbounds}
    \mch(r) \leq \min \left\{ \max \left( f_0^{\ell}(\underbrace{1,\dots,1}_{r-1}) + 1, \ell+r \right) \Bigg| \; 1 \leq \ell \leq H_R(r)-r \right\}.
\end{equation}
\end{theorem}

\begin{proof}
Consider $\tau_{1,\dots,r-1} \subset \Pb_{K_n}^{n-1}$. For any $r < n$, we have $\dim \tau_{1,\dots,r-1} \geq 1$, so showing that \[\tau_{1,\dots,r-1}(K_n^{(\ell)}) \subseteq \tau_{1,\dots,r-1}\] is dense means we can find $K_n^{(\ell)}$-points of $\tau_{1,\dots,r-1}$ distinct from $[1:0:\cdots:0]$. Finding such a point yields a valid Tschirnhaus transformation which eliminates $r$ variables at the cost of passing to an extension $L/K_n$ of resolvent degree at most $\ell$. Since $\tau_{1,\dots,r-1}$ is an intersection of type $(\underbrace{1,\dots,1}_{r-1})$ in $\Pb_{K_n}^{n-1}$, setting $N \coloneqq f_0^{\ell}(\underbrace{1,\dots,1}_{r-1}) + 1$ gives
\[ \RD(n) \leq \max( n-r, \ell ) \]
for any $n \geq N$. To prove the theorem, we consider this bound at $n = \max(N,\ell+r)$ to show that $\mch(r) \leq \max(N,\ell+r)$. We proceed in two cases. 

Firstly, if $N \geq \ell+r$, then at $n = \max(N,\ell+r) = N$ we have
\[ \RD(N) \leq \max(N-r,\ell) = N-r, \]
which implies $N - \RD(N) \geq r$, and therefore $\mch(r) \leq N$. 
On the other hand, if $N < \ell+r$, then we evaluate our bound at $n = \ell+r$ to see that
\[ \RD(\ell+r) \leq \max(\ell+r-r,\ell) = \ell. \]
This implies $(\ell+r) - \RD(\ell+r) \geq (\ell+r) - \ell = r$, and therefore $\mch(r) \leq \ell+r$. 
In both cases, we conclude $\mch(r) \leq \max(N, \ell+r)$. We then minimize over all valid $\ell \leq H_R(r)-r$ to obtain the sharpest possible bound.
\end{proof}

\begin{example}
Brauer's bound \eqref{eqn:brauer} amounts to carrying out the obliteration algorithm in one step. We show this via induction on $r$, where $\mch(4) = 5 \leq 3!+1$ is Bring's reduction of the quintic. Indeed, if we suppose that $\mch(r) \leq (r-1)!+1$ holds for some $r \geq 4$, then we have
\[ f^{r!-r}_0(\underbrace{1,\dots,1}_{r}) \leq f_r(0) = r, \]
via Lemma \ref{lemma:pts_from_planes}, since $r! \geq (r-1)!+1$ and thus $\RD(r!) \leq r! - r$ by assumption. Hence, Theorem \ref{thm:sharpening_bound_on_RD} gives
\[ \mch(r+1) \leq \max(r+1,(r!-r)+(r+1)) = r! + 1. \]
\end{example}

As outlined in Section \ref{sec:known_bounds}, our approach via Theorem \ref{thm:sharpening_bound_on_RD} is most effective at sharpening bounds in the ``intermediate range.'' Because our recursive methods rely on generic bounds, they fail to incorporate fine enumerative phenomena. Consequently, for very small values of $r$, our bounds can underperform classical constructions:

\begin{example}\label{example:not_always_best}
We saw in Example \ref{example:computations_for_nonic} that, for any $\ell < \RD(12)$, the best bound available to the methods in this paper for $f_0^\ell(1,1,1,1)$ is $9$; indeed, the same bound holds for $\mce = \Sol$ using the methods in \cite{GGW2025}. Taking $\ell = 1$, Theorem \ref{thm:sharpening_bound_on_RD} yields $\mch(5) \leq \max(9+1,1+5) = 10$, and so we can eliminate $5$ variables from the generic degree $n$ polynomial starting at $n = 10$. On the other hand, Hilbert's sketch \cite{Hilbert1927, Dixmier1993} uses the $27$ lines on a smooth cubic surface to show the strictly better bound $\mch(5) \leq 9$. 
\end{example}

Via explicit computation using the currently best known bounds on resolvent degree, we establish the new bounds on $\mch$ presented in Corollary \ref{cor:new_bounds_on_H_body}. Moreover, the new bounds for $\mch(8)$, $\mch(12)$, $\mch(13)$, $\mch(21)$ and $\mch(22)$ extend the bounds of the form $\frac{(r-1)!}{d!}$ from Sutherland \cite{Sutherland2022} given in terms of obtaining $(r-d+1)$-planes in $\tau_{1,\dots,d} \setminus \{ (1:0:\cdots:0) \}$, defined over extensions of sufficiently low resolvent degree:

\begin{corollary}[New Bounds on Resolvent Degree]\label{cor:new_bounds_on_H_body} The following bounds hold:
\begin{align*}
\mch(7) &\leq 74, & \mch(8) &\leq \tfrac{7!}{4!}+1, & \mch(11) &\leq 59050, \\
\mch(12) &\leq \tfrac{11!}{5!}+1, & \mch(13) &\leq \tfrac{12!}{5!}+1, & \mch(20) &\leq 227214539745187, \\
\mch(21) &\leq \tfrac{20!}{6!}+1, & & \text{and} & \mch(22) &\leq \tfrac{21!}{6!}+1.
\end{align*}
\end{corollary}

\begin{proof}
We proceed by casework for the given values of $r$. By Theorem~\ref{thm:sharpening_bound_on_RD}, we reduce this to a question of bounding $f_0^\ell$ through iteratively applying Corollary \ref{cor:obliteration} of Section \ref{sec:special_points_on_variety_of_j-planes}.

For the case $r=7$, we can use the bound $\mch(6)\le 21$ established in \cite{Sutherland2021} to compute the bound at $\ell=42$ using our greedy obliteration strategy \eqref{eqn:obliteration_strategy}:
\begin{gather*}
    f_0^{42}(1,1,1,1,1,1) \leq f_0^{42}(10,6,3,1) + 2 \leq f_0^{42}(34,9,1) + 5 \leq f_0^{42}(59,5) + 10 \leq f_{0}^{42}(64) + 10 = 74.
\end{gather*}
In the above, our applications of Corollary \ref{cor:obliteration} use the respective bounds
\[
    \RD(30) \leq 24, \qquad \RD(36) \leq 30, \qquad \RD(48) \leq 42, \qquad \text{and} \qquad \RD(32) \leq 26,
\]
which are all implied by $\mch(6)\le 21$. This gives $\mch(7) \le \max\{74+1, 42+7\} = 75$.

However, we can obtain a better bound at $\ell = 58$ via a different strategy. During the calculation of $f_0^{58}(1,1,1,1,1,1)$, the inequality reduces to bounding $f_0^{58}(34,9,1)$, similarly to the above, where our greedy strategy obliterates by $\ba = (0,4,1)$. If we instead take $\ba = (0,3,1)$, then we have
\begin{gather*}
    f_0^{58}(34,9,1) + 5 \leq f_0^{58}(58,6) + 9 \leq f_{0}^{58}(64) + 9 = 73.
\end{gather*}
We conclude $\mch(7) \leq \max\{73+1, 58+7\} = 74$ as claimed.

For $r=8$, we proceed similarly by using $\mch(7) \le 74$ to compute
\begin{gather*}
    f_0^{203}(1,1,1,1,1,1,1) \leq
    f_0^{203}(20,10,4,1) + 3 \leq
    f_0^{203}(70,14,1) + 7 \\
    \leq f_0^{203}(125,8) + 14 \leq
    f_0^{203}(133,1) + 21 \leq
    155.
\end{gather*}
The most costly obliteration is the first, where we use that $\RD(5\cdot 6\cdot 7) \leq 210-7 = \ell$.

The arguments for the remaining cases are similar, but cumbersome to reproduce in full detail; all use the greedy strategy \eqref{eqn:obliteration_strategy}. We leave these for the motivated reader (see Appendix \ref{appendix:computing_f}), with details regarding the computation of $f_0^\ell(1,\dots,1)$ at a minimizing value of $\ell$ for Theorem \ref{thm:new_bounds_on_H} shown in Table \ref{table:improvements-of-rd}.
\end{proof}

\begin{table}[htbp]
\centering
    \begin{tabular}{ccc}
        \toprule
        $r$ & $\ell$ & Bound on $f^{\ell}_0(1,\dots,1)$ \\
        \midrule
        7 & 58 & 73 \\
        8 & 203 & 155 \\
        11 & 59039 & 57367 \\
        12 & 332629 & 166777 \\
        13 & 3991668 & 400984 \\
        20 & 205891132094630 & 227214539745186 \\
        21 & 3379030566911980 & 625790049145201 \\
        22 & 70959641905151979 & 1603978712474280. \\
        \bottomrule
    \end{tabular}
\caption{Summary of calculations which improve bounds on $\mch(r)$. The tuple $(1,\dots,1)$ is understood to repeat $r-1$ times.}
\label{table:improvements-of-rd}
\end{table}

\begin{remark}\label{remark:hyperplanes_as_irrationalities}
Computing bounds on $N \coloneqq f^{\ell}_0(\overbrace{1,\dots,1}^{r-1})$ mirrors a process for writing down a formula for the generic degree $n \geq N$ polynomial in $n-r$ variables. Indeed, obliteration by $\ba$ is tantamount to introducing an accessory irrationality of degree $2^{a_2}\cdots n^{a_n}$. The computation for $\mch(7)$, for example, describes a tower of field extensions which makes sense for all $n \geq 75$:
\[ L_0 \coloneqq K_n \xhookrightarrow{\RD \leq 26} L_1 \xhookrightarrow{\RD \leq 42} L_2 \xhookrightarrow{\RD \leq 30} L_3 \xhookrightarrow{\RD \leq 24} L_4 \xhookrightarrow{\RD \leq n-7} L_5, \] 
where the extension $L_5/L_4$ corresponds to solving the reduced polynomial problem \eqref{eqn:reduced_poly}. The length of these formulas thus relate to the number of times Corollary \ref{cor:obliteration} is applied, which grows rapidly in $r$: e.g., Example \ref{example:many_accessory_irrationalities} gives a formula of length 1029977 for $\mch(22)$.
\end{remark}

In summary: Our methods match or improve upon the best estimates available in the literature for $6 \leq r \leq 33$ and $38 \leq r \leq 42$. Further, using the obliteration lemmas established in \cite{GGW2025}, the calculation $f^{\Sol}_0(1,1,1) \leq 4$ confirms the well-established result that the quintic can be reduced to a single-variable problem via a solvable Tschirnhaus transformation, i.e., $\mch(4) = 5$. However, the sharpest bounds from the obliteration method described in this paper give worse estimates than those in \cite{Sutherland2021} for $34 \leq r \leq 37$. As discussed in Appendix \ref{appendix:computing_f}, evaluating a bound for $f_0^\ell(\bm)$ becomes prohibitively expensive as $\ell$ approaches the minimizing value of \eqref{eqn:sharpenbounds}, because the number of necessary obliterations blows up so dramatically. As such, we could not determine if the methods in this paper \emph{improve} known bounds for $r \geq 40$ due to excessive runtimes.

\section{Bounds on \texorpdfstring{$\RD(G)$}{RD(G)} for finite groups}\label{sec:bounds_for_sporadics}

The introduction of $f_0^\ell(\bm)$, together with the accompanying computational lemmas in terms of polars, provides a uniform framework for the application of special points in \cite{GGSW2024} towards the resolvent degree of finite simple groups. 

\subsection{Torsors and \texorpdfstring{$\mce$}{E}-versality}\label{subsec:e-versality}

The argument in Section \ref{subsec:sporadic} follows \cite[Section 4]{GGSW2024}, which uses the language of torsors for algebraic groups over $k$. We retain this language for consistency with the surrounding literature, although our interest is the case where $G$ is an abstract finite group, viewed as a constant group scheme over $k$. We will assume $G$ is finite for the remainder of the paper, and refer the reader to \cite[Section~3]{GGSW2024} for additional details.

\begin{definition}[{\cite[Definition~2.1]{GGSW2024}}]
Let $X$ be a quasi-projective $G$-variety over $k$. A \emph{right $G$-torsor} over $X$ is a flat morphism $T \to X$ of $k$-schemes equipped with a right $G$-action $\sigma : G \times T \to T$ such that the natural map $G \times T \to T \times_X T$ given by $(g,t) \mapsto (\sigma(g,t),t)$ is an isomorphism. 
\end{definition}

We are primarily concerned with the case where the base scheme is $\Spec(K)$, with $K$ a finitely generated field extension of $k$. In this setting, a $G$-torsor $T \to \Spec(K)$ corresponds to the spectrum of a $G$-Galois algebra $A/K$. We say the torsor \emph{splits} over an extension $L/K$ if the base change
\[
T \times_{\Spec(K)} \Spec(L)
\]
is isomorphic to the trivial $G$-torsor $G \times \Spec(L)$, i.e., $A \otimes_K L \cong L^{|G|}$ as $L$-algebras equipped with the natural $G$-action.

\medskip
For our approach in bounding $\RD_k(G)$, it is helpful to recall the original geometric definition due to Farb--Wolfson \cite[Definition 3.1]{FW2019}:
\begin{equation}\label{eqn:def_rd_g_branched}
    \RD_k(G) = \sup\left\{ \RD_k(X \bra X/G) \st X \text{ is a primitive, faithful $G$-variety} \right\}.
\end{equation}
Note that, from our field-theoretic perspective, the branched cover $X \bra X/G$ induces a $G$-Galois extension of function fields $k(X)/k(X/G)$ such that
\[ 
    \RD_k(X \bra X/G) = \RD_k( k(X)/k(X/G) ). 
\]
For more on the equivalence of definitions, see \cite[Proposition~2.4]{FW2019}.

It is natural to ask for which $G$-varieties the supremum \eqref{eqn:def_rd_g_branched} is realized. A sufficient, though not necessary, condition is $X$ be \emph{$G$-versal}: that is, for every $G$-invariant Zariski open $U \subseteq X$ and every generically free $G$-variety $Y$, there is a $G$-equivariant rational map $Y \bra U$. For such varieties, we have $\RD_k(G) = \RD_k(X \bra X/G)$; moreover, the geometric version of the inequality arising from the trivial tower (Equation \ref{eqn:RD_ed_trdeg_bound} in Section \ref{subsec:ed_rd_special}) becomes
\begin{equation}\label{eqn:RD_g_dim_bound}
    \RD_k(G) \leq \dim X
\end{equation}
whenever $X$ is $G$-versal.

\begin{example}
The natural permutation action of $S_n$ on $\Ab^n$ is versal. Hence
\[
    \RD_k(n) = \RD\left(k(t_1,\dots,t_n)/k(t_1,\dots,t_n)^{S_n}\right),
\]
where the extension on the right is the splitting field of the generic polynomial of degree $n$. As we have seen in the many bounds on $\RD(n)$ discussed in Section \ref{sec:bounds_on_RD}, \eqref{eqn:RD_g_dim_bound} is not sharp here.
\end{example}

It was classically understood that versality is only a sufficient condition to realize the supremum of \eqref{eqn:def_rd_g_branched}. In particular, Klein used the icosahedral action $A_5 \acts \Pb^1_{\Cb}$ in his solution of the quintic \cite{Klein1884} to show that $\RD_\Cb(A_5) = 1$, while also showing that $\Pb^1_{\Cb}$ is not versal. Towards generalizing Klein's idea of solving an auxiliary quadratic, the notion of $\mce$-versality was introduced in \cite{FKW2023}. Following the framing of versality in terms of the arithmetic of rational points in \cite[Theorem~1.1]{DR2015}, $\mce$-versality was further developed in \cite{GGSW2024} to make explicit a technique involving level $\ell$ closures for producing bounds on $\RD_k(G)$. 

\begin{definition}[{\cite[Definition~3.2]{GGSW2024}}]
Let $X$ be an irreducible, generically free, quasi-projective $G$-variety over $k$, and let $\mce$ be a class of accessory irrationalities. We say that $X$ is:
\begin{itemize}
\item \emph{weakly $\mce$-versal for $G$} if for every $G$-torsor $T \to \Spec(K)$, there is an extension $K \into K'\in \mce(K)$ and a $G$-equivariant $k$-morphism
\[
    T \times_{\Spec(K)} \Spec(K') \to X;
\]
\item \emph{$\mce$-versal for $G$} if every $G$-invariant Zariski open $U\subset X$ is weakly $\mce$-versal.
\end{itemize}
\end{definition}

From this perspective, Klein showed that $\Pb_{\Cb}^1$ is \emph{solvably} versal, i.e., $\Sol$-versal, for $A_5$. 
Lastly, to state the arithmetic formulation of $\mce$-versality following \cite{DR2015}, we recall:

\begin{definition}[Twists]
    Given a $G$-torsor $T \to \Spec(K)$, $G$ acts on $T \times X$ diagonally and yields a $G$-torsor $T \times X \to \prescript{T}{}{X}$. We refer to $\prescript{T}{}{X}$ as the \emph{twist of $X$ by $T$}. 
\end{definition}

The core relationship between $K^{\mce}$-points and $\mce$-versality is as follows:

\begin{theorem}\cite[Theorem 3.9]{GGSW2024}\label{thm:twists}
    Let $G$ be a finite discrete group, and let $\mce$ be a class of accessory irrationalities which is saturated and closed under extensions. Let $X$ be an irreducible faithful $G$-variety over $k$. Then $X$ is $\mce$-versal if and only if for every finitely-generated $k$-field $K$ and $G$-torsor $T\to\Spec(K)$, we have that $\prescript{T}{}{X}(K^{\mce})\subset \prescript{T}{}{X}$ is dense. 
\end{theorem}

\subsection{New bounds for the sporadic groups}\label{subsec:sporadic}

In what follows, since $\RD_k(G) \leq \RD_{\Cb}(G)$ holds for any $k$ by \cite[Theorem~1.3]{Reichstein2025}, we restrict to $\RD(G) \coloneqq \RD_{\Cb}(G)$. We will describe an approach for bounding $\RD(G)$ for $G$ a finite simple group, then apply our procedure to the sporadics. Our interest in simple groups is due to \cite[Theorem~3.3 and Problem~3.5]{FW2019}, which establishes that $\RD(G)$ equals the max of $\RD(H)$ for all $H$ in the Jordan-H\"older decomposition of $G$. We write $\mu(G)$ for the degree of the minimal permutation representation of $G$.

\medskip
Given $G$ simple and an irreducible projective representation $\Pb(V) = \Pb^N$, we: 
\begin{itemize}
    \item identify $r$ algebraically independent invariant forms on $V$,
    \item obtain an $\RD_\Cb^{\leq \ell}$-versal $G$-subvariety $X \subseteq \Pb(V)$ for small $\ell$,
    \item bound $\RD(G)$ via $\RD(X \dashrightarrow X/G) \leq \dim X$.
\end{itemize}
The objective is to minimize $\ell$ while maximizing the number of forms $r$, using \cite[Proposition 3.3]{GGSW2024}:
\begin{equation}\label{eqn:rd_g_by_levels} 
    \RD(G) = \min_{\ell \geq 0} \left\{ \max\{\ell, \dim X\} \, \big| \, X \text{ is $\RD_{\Cb}^{\leq \ell}$-versal for } G \right\}. 
\end{equation}
Establishing $\RD_\Cb^{\leq \ell}$-versality uses Theorem~\ref{thm:twists}, which includes verifying generic freeness of the action $G \acts X$ and showing that, for every $G$-torsor $T \to \Spec K$ with $K$ finitely-generated over $\Cb$, each twist $\prescript{T}{}{X}$ has dense $K^{(\ell)}$-points. 

The approach in \cite{GGSW2024} to exhibit density uses a theorem of Merkurjev--Suslin together with a result of Reichstein (stated here for the case $k = \Cb$ and adapted to the context):

\begin{lemma}[{\cite[Lemma~14.4]{Reichstein2025}}] Let $K$ be finitely-generated over $\Cb$, with $X \subseteq \Pb^N$ a projective variety over $K$. If $\ell \geq \RD(\deg X)$, then $K^{(\ell)}$-points are dense in $X$. Moreover, if $Q \subseteq \Pb^N$ is a quadric of rank $r$ over $K$ with $\dim X \geq \floor{\frac{r+1}{2}}$, then $K^{(\ell)}$-points are dense in $X \cap Q$. 
\end{lemma}

In essence, the goal is to guarantee that $\RD(\deg X) \leq N-r$, with some minor improvements when one of the hypersurfaces defining $X$ is a quadric. With a more comprehensive approach to assess the density of $K^{(\ell)}$-points depending only on type, we can say more:

\begin{theorem}[Bounds on $\RD(G)$ via invariant theory]\label{thm:the_game} Let $G$ be a finite non-abelian simple group with a non-trivial irreducible projective representation $\Pb(V) = \Pb^N$.  Let $\mu(G)$ be the degree of the minimal permutation representation of $G$. Let $f_1,\dots,f_r$ be $G$-invariant forms on $\Pb(V)$ with respective degrees $d_1 \leq \dots \leq d_r$, such that each $f_i$ is algebraically independent from the lower degree invariants, and write $\bm$ for the corresponding type of $f_1,\ldots,f_r$. If $d_1 \cdots d_r < \mu(G)$ and $f_0^{N-r}(\bm) \leq N$, then 
\[ \RD(G) \leq N - r. \]
\end{theorem}
\begin{proof}
Our goal is to exhibit the $\RD_{\Cb}^{\leq N-r}$-versality of $X \coloneqq \Vb(f_1,\dots,f_r) \subseteq \Pb^N$, from which the result follows. First, because $G$ is simple and $\Pb(V)$ is a non-trivial projective representation, the $G$-action is faithful. Since $X$ has at most $d_1\dots d_r$ components which are permuted by $G$, we see that our assumption on $\mu(G)$ forces $X$ to be irreducible. Thus by \cite[Lemma 2.15]{GGSW2024} we conclude that $X$ is generically free. 

Moreover, the assumption on $f_0^{N-r}(\bm) \leq N$ means that $X$ has dense $K^{(N-r)}$-points; we need to confirm this holds for arbitrary twists $\prescript{T}{}{X}$ as in Theorem \ref{thm:twists}. Indeed, for any $G$-torsor $T \to \Spec K$ with $K$ finitely-generated over $\Cb$, the $G$-equivariant closed immersion $X \into \Pb(V)$ naturally induces the closed immersion $\prescript{T}{}{X} \into \prescript{T}{}{\Pb(V)}$. Since $\prescript{T}{}{\Pb(V)}$ is a Severi--Brauer variety which splits over $K^{\Sol}$ by the Merkurjev--Suslin Theorem \cite{MerkurjevSuslin1983}, it follows that $\prescript{T}{}{X}$ is an intersection of hypersurfaces in $\Pb^N_{K^{(N-r)}}$ of type $\bm$, and thus has dense $K^{(N-r)}$-points by assumption. We conclude the desired bound by \eqref{eqn:rd_g_by_levels}.
\end{proof}

We note that the premise $d_1 \cdots d_r < \mu(G)$ is only used to guarantee the action of $G$ on $\Vb(f_1,\dots,f_r)$ is generically free. There are certainly cases where a more careful analysis, i.e., not simply in terms of degrees of the underlying forms, could permit tighter bounds on $\RD(G)$. We also emphasize that the current best bounds on $\RD(d_1 \cdots d_r) + r$ are, in general, significantly larger than the corresponding bounds on $f_0^{N-r}(\bm)$, even when compared with the coarse polynomial bounds established in Section \ref{subsec:bounds_on_f}.

\begin{corollary}[New Bounds on $\RD(G)$ for the Sporadic Groups]\label{cor:sporadic_bounds_body} 
The following bounds hold:

\[ \begin{aligned}
\RD(\MCL) & \leq 18, \qquad & \qquad \RD(\RU) & \leq 25, \\
\RD(\HE) & \leq 47, & \RD(\FIWH) & \leq 775, \\
\RD(\FIWF) & \leq 778, & \RD(\B) & \leq 4364, \\
\RD(\MO) & \leq 196872.
\end{aligned} \]
\end{corollary}

\begin{proof}
The argument follows the analysis of \cite[Section 4]{GGSW2024}, employing additional invariants of select projective representations for the given sporadic groups using our more effective method for verifying the density of $K^{(d)}$-points. These invariants are summarized in Table \ref{table:sporadic-bounds}. The new bounds are obtained via the straightforward, though somewhat computationally intensive, calculations:
\[
\begin{array}{r@{}l@{}c@{}r@{}l}
    f_0^{18}(0,1,0,0,1,1)\ & \leq 21           && \RD(\MCL)\ & \leq 18, \\
    f_0^{25}(0,0,0,1,0,0,0,1)\ & \leq 5        && \RD(\RU)\ & \leq 25, \\
    f_0^{47}(0,0,1,1,1)\ & \leq 8              && \RD(\HE)\ & \leq 47, \\
    f_0^{775}(0,1,1,1,2,1)\ & \leq 24          &\qquad \Longrightarrow \qquad\quad& \RD(\FIWH)\ & \leq 775, \\
    f_0^{778}(0,0,1,0,0,2,0,0,1)\ & \leq 13    && \RD(\FIWF)\ & \leq 778, \\
    f_0^{4364}(0,1,0,1,0,1,0,3)\ & \leq 108    && \RD(\B)\ & \leq 4364, \\
    f_0^{196872}(0,1,1,1,1,3,3)\ & \leq 168825 && \RD(\MO)\ & \leq 196872.
\end{array}
\]
These are computed by applying Corollary \ref{cor:obliteration} according to our greedy strategy \eqref{eqn:obliteration_strategy}.
\end{proof}

\begin{table}[htbp]
\begin{tabular}{cccccc}
\toprule
$G$ & $\RD(\mu(G))$ bound & $\dim \Pb(V)$ & Invariants in \cite{GGSW2024} & New invariants & Thm. \ref{thm:the_game} bound \\
\midrule
\rowcolor{c6} $\JW$ &93 & 5 &None & None & 5\\
\rowcolor{c1} $\MOO$ &6 & 9 &2, 3 & None & 7 \\
\rowcolor{c1} $\MOW$ &7 & 10 &2, 3 & None & 8 \\
\rowcolor{c6} $\MWW$ &16 &9 &4 & None & 8 \\
\rowcolor{c6} $\SUZ$ &1782 & 11 &12 & None & 10 \\
\rowcolor{c6} $\JH$ &78 &17 &6 & None & 16 \\
\rowcolor{c1} $\MWH$ &17 &21 &2, 3, 4 & None & 18 \\
\rowcolor{c6} $\MWF$ &18 &22 &2, 3, 4, 5 & None & 18 \\
\rowcolor{c6} $\HS$ &93 &21 &2, 4, 5 & None & 18 \\
\rowcolor{c5} $\MCL$ &267 &21 &2, 5 & 6 & 18 \\
\rowcolor{c6} $\COH$ &268 &22 &2, 6 & None & 20 \\
\rowcolor{c6} $\COW$ &2291 &22 &2, 8 & None & 20 \\
\rowcolor{c6} $\COO$ &98269 &23 &2, 12 & None & 21 \\
\rowcolor{c5} $\RU$ &4051 &27 &4 & 8 & 25 \\
\rowcolor{c5} $\HE$ &2049 &50 &3, 4 & 5 & 47 \\
\rowcolor{c6} $\JO$ &258 &55 &2, 3, 4, 4 & None & 51 \\
\rowcolor{c6} $\FIWW$ &3501 &77 &2, 6, 8 & None & 74 \\
\rowcolor{c6} $\HN$ &1139988 &132 &2, 6, 7 & None & 129 \\
\rowcolor{c6} $\THO$ &143126986 &247 &2, 8, 8 & None & 244 \\
\rowcolor{c6} $\ON$ &122749 &341 &6, 6, 6 & None & 338 \\
\rowcolor{c5} $\FIWH$ &31661 &781 &2, 3, 4, 5, 5 & 6 & 775 \\
\rowcolor{c5} $\FIWF$ &306925 &782 &3, 6, 6 & 9 & 778\\
\rowcolor{c6} $\JF$ &173067375 &1332 &4, 6, 6, 7 & None & 1328 \\
\rowcolor{c6} $\LY$ &8835143 &2479 &6, 6, 6, 6 & None & 2475 \\
\rowcolor{c5} $\B$ &13571954984 &4370 &2, 4, 6, 8, 8 & 8 & 4364 \\
\rowcolor{c5} $\MO$ & 97239461142009185977 &196882 &2, 3, 4, 5, 6, 6, 6, 7 &7, 7 & 196872 \\
\bottomrule
\end{tabular}
\caption{Bounds on $\RD$ for the sporadic groups and invariants for select representations.}
\label{table:sporadic-bounds}
\end{table}

We summarize and compare this process with the bounds in \cite{GGSW2024} in Table \ref{table:sporadic-bounds} using the best currently available bounds on $\RD(n)$, including for sporadic groups whose resolvent degree bounds are not improved by the methods of this paper. For a thorough accounting of the associated representations and derivation of invariant forms, see \cite[Figures 2 and 7]{GGSW2024}. We color the table according to the rules:

\begin{itemize}
    \item Blue if the best bound on $\RD(G)$ comes (strictly) from the minimal permutation representation of $G$, i.e., $\RD(\mu(G))$ as in \cite[Lemma 3.13]{FW2019};
    \item Green if the methods described here furnish a better bound on $\RD(G)$ than \cite{GGSW2024}, which happens if our bounds on $f_0^{N-r}(\bm)$ allow us to include additional invariant forms (these are recorded in the relevant column) when defining $X$; or
    \item Yellow if the bound is not improved from \cite{GGSW2024}.
\end{itemize}

\newpage
\appendix
\section{Computing bounds on \texorpdfstring{$f^\ell_j(\bm)$}{flj(m)}}\label{appendix:computing_f} %

Here we review methods for computing bounds on $f^\ell_j(\bm)$, especially for obtaining bounds on resolvent degree via Theorem \ref{thm:sharpening_bound_on_RD}, and compare these to approximations as in Lemma \ref{lemma:coarse_upper_bound}. While the obliteration lemmas provide a deterministic recursive path to a bound, the computational complexity of the recursion grows rapidly with $|\bm| = m_1 + \dots + m_n$ and $n$. Following the greedy strategy \eqref{eqn:obliteration_strategy} in Section \ref{subsec:bounds_on_f}, we implement a top-down approach that maximizes the degree of obliterated forms.

\begin{algorithm}
\caption{Computing the greedy obliteration vector $\ba \leq \bm$.}\label{code:obliteration}
\begin{algorithmic}[1]
\Require Level of versality $\ell$ and type $\bm = (m_1, \dots, m_n) \in \Mt$
\State $\ba \gets (0, \dots, 0)$ \Comment{Initialize the obliteration vector}
\State $P \gets 1$ \Comment{Accumulates the product of degrees}
\State $d \gets n$ \Comment{Start obliterating highest degree forms first}
\While{$d \geq 1$} 
    \If{$a_d < m_d$}
        \If{$R_{\text{bound}}(P \cdot d) \leq \ell$} \Comment{Check if we can obliterate another form}
            \State $a_d \gets a_d + 1$
            \State $P \gets P \cdot d$
        \Else
            \State \textbf{break} \Comment{Bound reached; prevents obliterating lower-order terms}
        \EndIf
    \Else
        \State $d \gets d - 1$ \Comment{Degree $d$ exhausted; proceed to next lower degree}
    \EndIf
\EndWhile
\State \Return $\ba, \bm - \ba$ \Comment{Returns the obliteration and the remainder}
\end{algorithmic}
\end{algorithm}

We also include an implementation of the best known bound on $\RD$ via its Hamilton function (Algorithm \ref{code:RD} and Table \ref{table:H_bound}). Note that the values appearing in this Hamilton function incorporate our improvements in Theorem \ref{thm:new_bounds_on_H}. Bounds for $34 \leq r \leq 55$ are derived from \cite{Sutherland2021}\footnote{We alert the reader to a typo in the definition of $G$ on page 22, which should read ``$4 \leq d < m-1$''}. 

\begin{table}[htpb]
\centering
\renewcommand{\arraystretch}{1.2}
\begin{tabular}{cc}
\toprule
$r$ & Bound on $\mch(r)$ \\
\midrule
$1 \leq r \leq 4$ & $r + 1$ \\
$r = 5$ & $9$ \\
$r = 6$ & $21 = \frac{5!}{3!}+1$ \\
$r = 7$ & $75$ \\
$8 \leq r \leq 11$ & $\frac{(r-1)!}{4!} + 1$ \\
$12 \leq r \leq 19$ & $\frac{(r-1)!}{5!} + 1$ \\
$r = 20$ & $227214539745187$ \\
$21 \leq r \leq 33$ & $\frac{(r-1)!}{6!} + 1$ \\
$r = 34$ & $2475934708812781843231486891102123$ \\
$35 \leq r \leq 43$ & $\frac{(r-1)!}{7!} + 1$ \\
$r = 44$ & \quad $8559276927975810009082900078329761155025671771554$ \\
$45 \leq r \leq 55$ & $\frac{(r-1)!}{8!} + 1$ \\
\bottomrule
\end{tabular}
\caption{Explicit values of the Hamilton function for the best known bound on resolvent degree.}
\label{table:H_bound}
\end{table}

\begin{algorithm}
\caption{Implementation of $R_{\text{bound}}(n)$ using binary search on $H_{\text{bound}}$.}\label{code:RD}
\begin{algorithmic}[1]
\Require Array $H_{\text{bound}}$ precomputed for desired values of $r$ (e.g., by Table \ref{table:H_bound}).
\State $r_{\max} \gets \max(\{r \mid H_{\text{bound}}[r] \text{ is defined}\})$
\vspace{0.5em}
\Function{$R_{\text{bound}}$}{$n$}
    \If{$n < H_{\text{bound}}[1]$}
        \State \Return $n$
    \EndIf
    \If{$n \geq H_{\text{bound}}[r_{\max}]$}
        \State \Return $n - r_{\max}$
    \EndIf
    \State $\texttt{low} \gets 1$
    \State $\texttt{high} \gets r_{\max}$
    \While{$\texttt{low} \leq \texttt{high}$} \Comment{Search for $r$ with $H[r] \leq n < H[r+1]$}
        \State $r \gets \lfloor (\texttt{low} + \texttt{high}) / 2 \rfloor$
        \If{$H_{\text{bound}}[r] \leq n$ \textbf{and} $H_{\text{bound}}[r+1] > n$}
            \State \Return $n - r$
        \ElsIf{$H_{\text{bound}}[r] \leq n$}
            \State $\texttt{low} \gets r + 1$
        \Else
            \State $\texttt{high} \gets r - 1$
        \EndIf
    \EndWhile
\EndFunction
\end{algorithmic}
\end{algorithm}

In principle, one can attempt to sharpen a bound on $\mch(r)$ using Theorem \ref{thm:new_bounds_on_H} by performing a binary search for the optimal level $\ell$ in the range $[1, H_{\text{bound}}(r)-r]$. Since $f^\ell_0(\bm)$ is non-increasing in $\ell$, the objective is to find the minimum $\ell$ such that $f^\ell_0(1,\dots,1) + 1 \leq \ell + r$. 

In practice, this approach becomes computationally intractable when $\ell$ falls below a certain threshold. A smaller level $\ell$ restricts the size of the obliteration vector $\ba$ at each step, which compounds dramatically after successive applications of endomorphism $\bm \to \bm^j$. These factors significantly increase the number of iterations required to reach the base case. For the larger values of $r$ considered in this paper, choosing $\ell$ even slightly below the optimal value results in a chain of recursive calls that exceeds standard memory and time constraints. 

The following example compares the approximation in terms of obliteration by lines in Lemma \ref{lemma:coarse_upper_bound} with our greedy obliteration strategy, illustrating the sensitivity of the latter process to $\ell$ as $n$ grows. For the sake of consistency with the remainder of the paper, these calculations use the \emph{previously} best available bounds:
\begin{align*}
\mch(7) &\leq 109, & \mch(8) &\leq 325, & \mch(11) &\leq \tfrac{10!}{4!} + 1, \\
\mch(12) &\leq \tfrac{11!}{4!} + 1, & \mch(13) &\leq 5250198, & \mch(20) &\leq \tfrac{19!}{5!} + 1, \\
\mch(21) &\leq \tfrac{20!}{5!} + 1, & & \text{and} & \mch(22) &\leq 381918437071508900.
\end{align*}

\begin{example}\label{example:coarse_bound_is_coarse}
From Lemma \ref{lemma:coarse_upper_bound}, we have the coarse bound $f_0^\ell(0,m_2) \leq \tbinom{m_2+1}{2}$. Using $\ell = \RD(2^{100}) \leq 2^{100}-31$ (\cite{Sutherland2022}), we compare to the best bounds available in this paper.
\begin{center}
    \begin{tabular}{r|cccccc}
        \toprule
        $m_2$ & 100 & 101 & 120 & 360 & 1080 & 4320 \\
        \midrule
        Best bound on $f_0^\ell(0, m_2)$ & 100 & 201 & 2120 & 48360 & 531080 & 9120320 \\
        Bound via Lemma \ref{lemma:coarse_upper_bound} & 5050 & 5151 & 7260 & 64980 & 583740 & 9333360. \\
        \bottomrule
    \end{tabular}
\end{center}
As another point of comparison, consider the greedy obliteration strategy with our coarse bound on the intersection of $m_4$ quartics with the same $\ell = \RD(4^{50})$:
\begin{center}
    \begin{tabular}{r|cccc}
        \toprule
        $m_4$ & 50 & 51 & 60 & 120 \\
        \midrule
        Greedy bound on $f_0^\ell(0,0,0,m_4)$ & $50$ & $\sim 8.39 \cdot 10^{5}$ & $\sim 8.11 \cdot 10^{9}$ & $\sim 5.37 \cdot 10^{13}$ \\
        Bound via Lemma \ref{lemma:coarse_upper_bound} & $\sim 3.14 \cdot 10^{11}$ & $\sim 3.67 \cdot 10^{11}$ & $\sim 1.34 \cdot 10^{12}$ & $\sim 3.40 \cdot 10^{14}$. \\ 
        \bottomrule
    \end{tabular}
\end{center}
\end{example}

\subsection{Elementary speed-ups}\label{appendix:elementary_speed-ups}

The sensitivity of the obliteration method to $\ell$ can make searching for the minimizing value of Theorem \ref{thm:sharpening_bound_on_RD} computationally expensive.
While implementing a bailout threshold can marginally address this issue, we outline two more effective strategies. 

Firstly, as shown in Example \ref{example:many_accessory_irrationalities}, the number of calculations needed to obtain estimates grows dramatically with $n$ and becomes concentrated in the lower degree terms. Given the high volume of maximally quadratic obliterations that must be evaluated in these recursions, we provide a closed-form formula for the direct computation of $f_j^\ell(m_1, m_2)$ to bypass several layers of recursive calls.

\begin{proposition}\label{prop:l-points_on_quadrics}
Fix $m_1, m_n \in \Nb$. For $b, \ell \in \Nb_{\geq 1}$ with $\RD(2^b) \leq \ell$, write $m_2 = q b + r$ with $0\le r<b$ and $q \in \Nb$ as in the division algorithm. Then
\[ f_0^\ell(m_1,m_2) \leq m_1 + \binom{q}{2} b^2 + q b (r+1) + r. \] 

\end{proposition}
\begin{proof}
We induct on $q$. The base case $q=0$ (where $m_2 = r < b$) follows from a single use of Lemma \ref{lemma:pts_from_planes}:
\[ f^\ell_0(m_1, r) \leq f^\ell_{m_1+r}(0) = m_1 + r. \]
Supposing the formula holds for $m_2 = qb + r$, we consider $m_2' = (q+1)b + r$. By Corollary \ref{cor:obliteration} and the bound $\RD(2^b) \leq \ell$, we can obliterate $b$ quadrics to obtain:
\[ f^\ell_0(m_1, m_2') \leq f^\ell_b(m_1, m_2) \leq f^\ell_0(m_1 + b m_2, m_2) + b. \]
Substituting and applying the inductive hypothesis gives:
\[
\begin{split}
    f^\ell_0(m_1, m_2') & \leq (m_1 + qb^2 + br) + \tbinom{q}{2} b^2 + q b (r+1) + r + b \\
    & = m_1 + \left( \tbinom{q}{2} b^2 + q b^2 \right) + \left( qb(r+1) + b(r+1) \right) + r \\
    & = m_1 + \tbinom{q+1}{2} b^2 + (q+1) b(r+1) + r. \\
\end{split}
\]
\end{proof}

\noindent While this speed-up sufficed for our purposes, such results can be extended arbitrarily to evaluation of $f_j^\ell(m_1,\dots,m_n)$ given sufficient interest and knowledge of the growth of $\RD(d^{b_d})$ for $2 \leq d \leq n$.

Secondly, rather than implementing a naive search to find an optimal $\ell$, we use the following observation (due to the Hamilton condition, Definition \ref{def:bound_on_RD}\ref{def:hamilton_condition}): as $\ell$ decreases, the sequence of obliterations for a given $\bm$ only changes when the bound $R_{\text{bound}}$ prevents us from extracting the same factors. Specifically, by keeping track of the maximum value of $P = 2^{a_1} \cdots n^{a_n}$ (the product of degrees) encountered in Algorithm \ref{code:obliteration} during the evaluation of $f^\ell_0(\bm)$, we can identify the next ``critical'' value of $\ell$: Our next guess in the optimization process of Theorem \ref{thm:sharpening_bound_on_RD} is $\ell = R_{\text{bound}}(P) - 1$, i.e., the first level where we can no longer obliterate the $\ba$ corresponding to the maximum $P$.

We conclude by illustrating these optimizations with two examples:

\begin{example}
Writing $R$ for the best bounds on $\RD$ before this paper, we consider $f^\ell_0(\overbrace{1,\dots,1}^{20})$ in sharpening $\mch(21)$. Table \ref{table:improvements-of-h-21} includes the values of $\ell$ tested, the largest product extracted over the course of each calculation (which is used to determine the next value of $\ell$), and a count of the number of inequalities needed to reach either the base case $f_j(0) = j$ (i.e. the length of the tower derived as in Remark \ref{remark:hyperplanes_as_irrationalities}) or the quadric case as in Proposition \ref{prop:l-points_on_quadrics}.

\begin{table}[htbp]
\centering
    \footnotesize
    \setlength{\tabcolsep}{10pt}
    \renewcommand{\arraystretch}{1.5}
    \begin{tabular}{ccccc}
        \toprule
        \small $\ell$ & \small $f^\ell_0(1,\dots,1) \leq$ & \small $\max(P)$ & \small Steps to $f_j(0)$ & \small Steps to $f_j(m_1,m_2)$ \\
        \midrule
        $\tfrac{20!}{5!}-20 $ & 52111719 & $\tfrac{20!}{5!}$ & 193 & 12 \\
        20274183401471979 & 580843710708473 & $2^{24} \cdot 3^{19}$ & 631321 & 247 \\
        19499511680335851 & 580843710708494 & $2^{54}$ & 631321 & 247 \\
        18014398509481963 & 604121162689288 & $2^{46} \cdot 3^5$ & 655989 & 250 \\
        17099604835172331 & 604121180066701 & $3^{34}$ & 655989 & 250 \\
        16677181699666548 & 604266478377039 & $2^{17} \cdot 3 \cdot 5^{15}$ & 656075 & 257 \\
        11999999999999979 & 608622926700042 & $2^{53}$ & 658436 & 257 \\
        8226356490141675 & 608622874490554 & $3^{33}$ & 671093 & 257 \\
        5559060566555502 & 608769515586248 & $2^{11} \cdot 3^{26}$ & 671181 & 265 \\
        5205741216417771 & 608769515586235 & $2^{52}$ & 671181 & 265 \\
        4503599627370475 & 625790031459181 & $2^{12} \cdot 3^{25}$ & 693836 & 267 \\
        3470494144278507 & 625790049145201 & $\frac{20!}{6!}$ & 693836 & 267 \\
        3379030566911979 & $\sim 2.795 \cdot 10^{28}$ & $2^{18} \cdot 3^{21}$ & 4636021784637 & 679683 \\
        \bottomrule
    \end{tabular}
\caption{Summary of calculations for improving the bound on $\mch(21)$. The tuple $(1,\dots,1)$ in the second column header is understood to repeat $20$ times.}
\label{table:improvements-of-h-21}
\end{table}

From the table, we can see how abruptly the calculation changes; moreover, the optimal value of $\ell$ occurs between the last two rows. The minimum of Theorem \ref{thm:sharpening_bound_on_RD} is achieved by $f_0^\ell(1,\dots,1)+1$ or by $\ell+r$, and our bound on $f^\ell_0(1,\dots,1)$ is constant for all $\ell$ in the range:
\[ R_{\text{bound}}(\tfrac{20!}{6!}) \leq \ell \leq 3470494144278507. \]
Since $f^\ell_0(1,\dots,1) + 1 = 625790049145202$ (as in the penultimate row) is smaller than any $\ell$ in this range, the optimal bound is obtained by choosing the smallest such $\ell$:
\[ \mch(21) \leq \ell + 21 = R_{\text{bound}}(\tfrac{20!}{6!}) + 21 = \left( \tfrac{20!}{6!} - 20 \right) + 21 = \tfrac{20!}{6!} + 1. \]
\end{example}

\subsubsection*{Data Availability}
The data and algorithms supporting the findings in this paper, together with additional references and supplementary information files, are available on GitHub \cite{GomezGonzales2025}. Further materials are available upon request.

\clearpage
\bibliography{sn-bibliography}

\end{document}

%% file: macros.tex
\newcommand{\catname}[1]{{\textsf{\textbf{#1}}}}

\newcommand{\Fieldsk}{\catname{Fields/k}}

\newcommand{\Cat}{\catname{Cat}}

\newcommand{\Ab}{\mathbb A}

\newcommand{\Cb}{\mathbb C}

\newcommand{\Nb}{\mathbb N}

\newcommand{\Pb}{\mathbb P}
\newcommand{\Qb}{\mathbb Q}

\newcommand{\Vb}{\mathbb V}

\newcommand{\Mt}{M_{\operatorname{type}}}
\newcommand{\ba}{\mathbf{a}}
\newcommand{\bb}{\mathbf{b}}
\newcommand{\be}{\mathbf{e}}
\newcommand{\bm}{\mathbf{m}}

\newcommand{\mcc}{\mathcal C}

\newcommand{\mce}{\mathcal E}
\newcommand{\mcf}{\mathcal F}

\newcommand{\mch}{\mathcal H}
\newcommand{\mci}{\mathcal I}

\newcommand{\lbm}{\left[ \begin{matrix}}
\newcommand{\rem}{\end{matrix} \right]}

\DeclareMathOperator{\PGL}{PGL}

\DeclareMathOperator{\grass}{Gr}

\newcommand{\st}{ \ | \ }

\newcommand{\floor}[1]{\left\lfloor {#1} \right\rfloor}
\newcommand{\acts}{\;\rotatebox[origin=c]{-90}{$\circlearrowright$}\;}

\newcommand{\into}{\hookrightarrow}

\newcommand{\bra}{\dashrightarrow}

\DeclareMathOperator{\Spec}{Spec}
\DeclareMathOperator{\red}{red}

\DeclareMathOperator{\Fin}{Fin}

\DeclareMathOperator{\Gal}{Gal}
\DeclareMathOperator{\trdeg}{tr.deg}
\DeclareMathOperator{\characteristic}{char}
\newcommand{\Sol}{\mathrm{Sol}}

\DeclareMathOperator{\ed}{ed}
\DeclareMathOperator{\RD}{RD}

\DeclareMathOperator{\MOO}{M_{11}}
\DeclareMathOperator{\MOW}{M_{12}}
\DeclareMathOperator{\MWW}{M_{22}}
\DeclareMathOperator{\MWH}{M_{23}}
\DeclareMathOperator{\MWF}{M_{24}}
\DeclareMathOperator{\COO}{Co_1}
\DeclareMathOperator{\COW}{Co_2}
\DeclareMathOperator{\COH}{Co_3}
\DeclareMathOperator{\JO}{J_1}
\DeclareMathOperator{\JW}{J_2}
\DeclareMathOperator{\JH}{J_3}
\DeclareMathOperator{\JF}{J_4}
\DeclareMathOperator{\SUZ}{Suz}
\DeclareMathOperator{\MCL}{McL}
\DeclareMathOperator{\HS}{HS}
\DeclareMathOperator{\FIWW}{Fi_{22}}
\DeclareMathOperator{\FIWH}{Fi_{23}}
\DeclareMathOperator{\FIWF}{Fi_{24}'}
\DeclareMathOperator{\THO}{Th}
\DeclareMathOperator{\HN}{HN}
\DeclareMathOperator{\HE}{He}
\DeclareMathOperator{\B}{B}
\DeclareMathOperator{\MO}{M}
\DeclareMathOperator{\ON}{O'N}
\DeclareMathOperator{\RU}{Ru}
\DeclareMathOperator{\LY}{Ly}
